\documentclass[aos,MSNbibl,nameyear,seceqn,dvips]{arximspdf}
\usepackage{dcolumn}

\doi{10.1214/12-AOS1005} %kopijuoti is PTS
\volume{40}
\issue{2}
\pubyear{2012}
\firstpage{1102}
\lastpage{1131}

\makeatletter
\newcolumntype{d}[1]{D{.}{.}{#1}}
\newtheorem{theorem}{Theorem}[section]
\newtheorem{lemma}[theorem]{Lemma}
\newtheorem{corollary}[theorem]{Corollary}

\newproclaim{assumption}{Assumption}[section]
\newproclaim{remark}{Remark}[section]
\newcommand{\bzero}{\mathbf{0}}

\def\bsuffix #1{#1}
\makeatother

\begin{document}
\begin{frontmatter}

\title{Bayesian empirical likelihood for quantile~regression\thanksref{T1}}
\runtitle{Bayesian empirical likelihood for quantile regression}
\thankstext{T1}{Supported in part by NSF Grants DMS-07-24752 and DMS-10-07396.}

\begin{aug}
\author[A]{\fnms{Yunwen} \snm{Yang}\ead[label=e1]{Yunwen.Yang@drexel.edu}}
\and
\author[B]{\fnms{Xuming} \snm{He}\corref{}\ead[label=e2]{xmh@umich.edu}}
\runauthor{Y. Yang and X. He}
\affiliation{Drexel University and University of Michigan}
\address[A]{Department of Epidemiology and Biostatistics\\
School of Public Health\\
Drexel University\\
1505 Race Street\\
Philadelphia, Pennsylvania 19102\\
USA\\
\printead{e1}} %adresu isvedimo komanda gale!
\address[B]{Department of Statistics\\
University of Michigan\\
1085 South University Avenue\\
Ann Arbor, Michigan 41809\\
USA\\
\printead{e2}}
\end{aug}

% HISTORY:
\received{\smonth{9} \syear{2011}}
\revised{\smonth{2} \syear{2012}}

% ABSTRACT
%
\begin{abstract}
Bayesian inference provides a flexible way of combining data with prior
information. However, quantile regression is not equipped with a~parametric likelihood, and therefore, Bayesian inference for quantile
regression demands careful investigation. This paper considers the
Bayesian empirical likelihood approach to quantile regression. Taking
the empirical likelihood into a Bayesian framework, we show that the
resultant \emph{posterior} from any fixed prior is asymptotically
normal; its mean shrinks toward the true parameter values, and its
variance approaches that of the maximum empirical likelihood estimator.
A~more interesting case can be made for the Bayesian empirical
likelihood when informative priors are used to explore commonality
across quantiles. Regression quantiles that are computed separately at
each percentile level tend to be highly variable in the data sparse
areas (e.g., high or low percentile levels). Through empirical
likelihood, the proposed method enables us to explore various forms of
commonality across quantiles for efficiency gains. By using an MCMC
algorithm in the computation, we avoid the daunting task of directly
maximizing empirical likelihood. The finite sample performance of the
proposed method is investigated empirically, where substantial
efficiency gains are demonstrated with informative priors on common
features across several percentile levels. A theoretical framework of
shrinking priors is used in the paper to better understand the power of
the proposed method.
\end{abstract}

% KEYWORDS
%


\end{frontmatter}

%s1 ###
\section{Introduction}

Quantile regression is a statistical methodology for the modeling and
inference of conditional quantile functions. Following
\citet{regquant}, we specify the $\tau$th conditional quantile function
of $Y \in\mathbb{R}$ given $X \in\mathbb{R}^{p+1} $ as
%e1 ###
%
\begin{equation}\label{quantregmd}
Q_{\tau}(Y|X)=X^\top\beta(\tau),
\end{equation}
where $\tau\in(0,1 )$, and $\beta(\tau)$ typically includes an
intercept. Quantile modeling of this type can be estimated for one or
several percentile levels; we refer the details on computation and
basic asymptotic theory to \citet{quantilereg}. Inferential methods for
quantile regression have been developed by a number of researchers,
including \citet{regrankscore}, \citet{pairbtstrap}, \citet{absdev} and
\citet{mcmbci}. The $\tau$-specific models allow for great flexibility,
as $\beta(\tau)$ for upper or lower quantiles can be distinct from
central trends, but the quantile estimates are highly variable in
data-sparse areas. Taking advantage of some commonality in the quantile
coefficients $\beta(\tau)$ across $\tau$ can provide a desirable
balance in the bias-variance tradeoff. In this article, we consider
using prior information on~$\beta(\tau)$ across several $\tau$
values. For example, a common slope assumption for $\tau$ near 1 can
improve the efficiency of high quantile estimation. Other forms of
informative priors on~$\beta(\tau)$ may achieve a similar goal.
Bayesian methods are a natural way of combining data with prior
information. The main difficulty in putting the Bayesian method to work
for quantile regression is that the model on $Q_{\tau}(Y|X)$ for one or
any small number of $\tau$ values does not specify a parametric
likelihood, which is needed in the Bayesian framework.

Several authors have attempted to use a working likelihood in the
Bayesian quantile regression framework. \citet{mixnonpara} and
\citet{mixnonpara2} used Dirichlet process mixture models.
\citet{bqrclu} assumed the error distributions to be an infinite
mixture of normals. \citet{extlatentmed} used an approximate method
based on the Jefferey's substitution likelihood for quantiles.
\citet{bayesdeqr}, \citet{qrlon} and \citet{bqrmixed}, among others,
chose (asymmetric) Laplace distributions as the working likelihood.
%model.
Those approaches, mostly tailored toward a specific percentile level of
$\tau$, use Markov chain Monte Carlo algorithms as a useful means of
computation. Work of these authors provided numerical evidence that a
Bayesian approach to quantile regression has merits.

In this article, we focus on estimating several quantiles together. To
do so, we use the empirical likelihood (EL), introduced by
\citet{empsigfunc}, to incorporate quantile regression into a
(pseudo-) Bayesian framework. Empirical likelihood makes it easy to
model several quantiles at the same time, allowing informative priors
on $\beta(\tau)$ across $\tau$ to be utilized. Statistical inference
based on empirical likelihood is known to enjoy good asymptotic
properties, especially if the EL is associated with moment restrictions
of sufficient smoothness. \citet{empnonsm} considered the EL with
nonsmooth estimating equations under a general setting. A~more
comprehensive review about empirical likelihood can be found in
\citet{OwenEL} and \citet{SongELreview}. Since the moment restrictions
for quantiles are placed on nonsmooth functions, some researchers,
including \citet{empsmq}, \citet{empsmqreg} and \citet{condemplik}
proposed using smoothed versions of the quantile estimating equations.
The smoothed EL is further extended to weakly dependent processes in
\citet{SongELdep} and censored data in \citet{RenELcensored}. We choose
to focus on the exact moment conditions for quantiles without the
complication of choosing a smoothing parameter. Those moment conditions
are also used in \citet{JudyELqr} and \citet{MiokEL} for clustered
data. In addition, we use a standard MCMC algorithm to explore the
posterior, to avoid the daunting task of directly maximizing the
empirical likelihood. In fact, the EL function given any proposed
parameters is relatively easy to compute, even though the
EL-maximization is notoriously difficult, even in modest dimensions.

The empirical likelihood is not a likelihood in the usual sense, so the
validity of the resultant posterior does not follow automatically from
the Bayes formula. \citet{bellazar} discussed the validity of inference
for the Bayesian empirical likelihood (BEL) approach based on earlier
work of \citet{covproperty}. \citet{betelest} and \citet{betelqr}
considered Bayesian exponentially tilted empirical likelihood (ETEL),
which can be viewed as a nonparametric Bayesian procedure with
noninformative priors on the space of distributions. \citet{betelqr}
further considered Bayesian ETEL in quantile regression. For the
inference of population means, \citet{secondcov} investigated the
asymptotic validity and accuracy of the Bayesian credible regions, and
furthermore, \citet{ByfrCI} showed that EL admits posterior based
inference with the frequentist asymptotic validity, but many of its
variants do not enjoy this property.
%To establish the asymptotic validity of the BEL for quantile
%regression, we need to work with discontinuous quantile estimating
%equations, so direct local expansions used in the EL literature cannot
%be used.
In this article, we establish the asymptotic distributions of the
posterior from the BEL approach for quantile regression, which enable
us to evaluate efficiency gains from informative priors.
\citet{lteestMCMC} discussed the asymptotic properties of the
quasi-posterior distributions defined as transformations of general
statistical criterion functions. In our work, we establish the
asymptotic distributions of the posterior from the BEL approach for
quantile regression, and are particularly interested in the interaction
of informative priors and empirical likelihood on the asymptotic
distribution of the posterior, which enables us to evaluate efficiency
gains from informative priors.
%Although finite-sample validity of the BEL posterior inference cannot
%be expected in our setting, we continue to use the term ``posterior"
%throughout the article.

Ideas similar to BEL have been used by other researchers.
\citet{baygmm} proposed the Bayesian generalized method of moments
(GMM), which can be adapted to quantile estimation. \citet{btstrp}
considered Bayesian bootstrap in quantile regression. Note that the GMM
estimators are also defined through moment restrictions, which allow
them to model multiple quantiles jointly. The GMM estimators, the
maximum empirical likelihood estimators (MELE) and some other EL-type
estimators generally have the same asymptotic distributions, but
possibly different higher order asymptotic properties; see \citet
{hordergel} and \citet{etelsche}. As discussed in \citet{hordergel},
the empirical likelihood approach has advantages over the GMM
estimators. Unlike GMM, the (asymptotic) bias of the MELE does not grow
with the number of moment restrictions. Furthermore, the efficiency of
the GMM estimator relies on a covariance matrix estimate for the
estimating equations, which could be ill-conditioned when estimating
multiple quantiles.

The recent development of Bayesian (conditional) density estimation
using mixture models enables nonparametric regression models on all
quantiles simultaneously; see \citet{Muller96}, \citet{Muller04},
\citet{Dunson07} and \citet{Dunson09}, among others. Theoretical
results about posterior consistency can be found in \citet{Pati10},
\citet{Norets11} and the references therein. In contrast, our proposed
BEL approach targets a small number of selected quantiles without the
need to model the entire conditional distributions. A novel part of our
work is its ability of employing informative priors to explore
commonality across quantiles for efficiency gains.

The rest of the paper is organized as follows. In Section
\ref{methodint}, we introduce the proposed BEL approach for quantile
regression, and discuss
model assumptions, method of computation and use of informative
priors. The asymptotic properties on the BEL posteriors
are provided in Section \ref{sec:asymprop} for both fixed and a class
of shrinking priors. The theoretical framework of shrinking priors
enables us to understand the efficiency gains of the BEL approach over
traditional methods. Section \ref{sec:simulation} demonstrates the
finite sample performance of the BEL approach
through Monte Carlo simulations with a focus on frequentist properties
of BEL posterior intervals, and efficiency gains from informative
priors. In Section \ref{sec:application}, we use a real data example
to show that the BEL approach
can be used as a useful statistical
downscaling method for the projection of high quantiles of temperature
from large scale climate models to a local scale. Some concluding
remarks are given in Section \ref{sec:discussion}.
The technical details to support the theorems in Section \ref{sec:asymprop}
are provided in the \hyperref[sec-proof]{Appendix}.

%s2 ###
\section{Bayesian empirical likelihood for quantile regression} \label
{methodint}

In this section we introduce the Bayesian empirical likelihood approach
for quantile regression. We begin with notation and definitions of the
underlying models and moment restrictions. Let $D=
\{(X_i,Y_i),i=1,\ldots,n \}$ be a random sample from the following
quantile regression model:
%e2 ###
%
\begin{equation}\label{exp-linearqrmodel}
Q_{\tau}(Y|X) =X^\top\beta_0(\tau),
\end{equation}
where $X \in\mathbb{R}^{p+1}$ is composed of an intercept term and $p$
covariates.
%The design points are either fixed or random.
We assume that the distribution of the $p$ covariates, $G_X$, has a
bounded support~$\mathcal{X}$. If the design points are nonstochastic,
the basic conclusions we obtain in this paper hold under appropriate
conditions on the design sequence, but we focus on the case of random
designs for simplicity. The unknown function~$\beta_0(\tau)$, if
specified over all $\tau\in(0,1)$, describes the entire conditional
distribution of $Y$ given $X$, which is denoted as $F_X$ in the rest of
the paper. We consider the problem of estimating $k$ quantiles at
$\tau_1<\tau_2<\cdots<\tau_k$, and let
$\zeta_{0}=(\beta_0(\tau_1),\ldots,\beta_0(\tau_k))$ be the true parameter
of interest in $\mathbb{R}^{k(p+1)}$. In most applications, $k$ is a
small integer. To estimate $\zeta_{0}$, we use $k(p+1)$ dimensional
estimating functions $m(X,Y,\zeta)$, where
$\zeta=(\beta(\tau_1),\ldots,\beta(\tau_k))$ and the components of $m$ are
%e3 ###
%
\begin{equation}\label{esteqael}
m_{dk + j}( X,Y,\zeta) =
\psi_{\tau_{d+1}}\bigl(Y-X^\top\beta(\tau_{d+1})\bigr)X_j
\end{equation}
for $d=0,1,\ldots,k-1$, $j=0,1,\ldots,p$, with
\[
\psi_{\tau}(u)=
\cases{
 1_{\{u<0\}}-\tau, &\quad $u \neq0$, \cr
0, &\quad $u=0$
}
\]
being the quantile score function, where $1_{\{A\}}$ is an indicator
function on the set~$A$. We hasten to add that $\zeta$ may contain
fewer than $k(p+1)$ unknown parameters when some common parameters are
present in $\beta(\tau)$ at different quantile levels. In such cases,
the number of moment restrictions exceeds the number of unknown
parameters. As shown in \citet{empgee} for smooth estimating functions,
the maximum empirical likelihood estimator attains the optimal
asymptotic efficiency subject to those moment conditions. We expect the
same for quantile estimating functions.

For any proposed $\zeta$, its profile empirical likelihood ratio is
given by
%e4 ###
%
\begin{equation}\label{empratio}\qquad
\mathcal{R}(\zeta) = \max \Biggl\{ \prod_{i=1}^n
(n\omega_i) \bigg| \sum_{i=1}^n \omega_i m(X_i, Y_i, \zeta)=0,
\omega_i \geq0, \sum_{i=1}^n \omega_i =1 \Biggr\}.
\end{equation}
By a standard Lagrange multiplier argument, we have
\[
\mathcal{R}(\zeta)= \prod_{i=1}^n \{ n\omega_i(\zeta) \},
\]
where the weights $\omega_i(\zeta) = [n\{1+\lambda_n(\zeta)^\top
m(X_i,Y_i,\zeta)\}]^{-1}$, and the Lagrange multiplier
$\lambda_n(\zeta)$ satisfies the following equation:
\[
\sum_{i=1}^n\frac{m(X_i,Y_i,\zeta)}{1+\lambda_n(\zeta)^\top
m(X_i,Y_i,\zeta)}=0.
\]
As discussed in \citet{extlam} and \citet{empgee}, the existence and
uniqueness of $\lambda_n(\zeta)$ are guaranteed when the following two
conditions are satisfied:
\begin{longlist}
\item[(C1)] The vector $\underline{0}\in\mathbb{R}^{k(p+1)} $ is
within the convex hull of $ \{m(X_i,Y_i,\zeta),\allowbreak i=1,\ldots,n \}$.
\item[(C2)] The matrix $\sum_{i=1}^{n}
\{m(X_i,Y_i,\zeta)m(X_i,Y_i,\zeta)^\top\}$ is positive\vadjust{\goodbreak}
definite.
\end{longlist}
The first condition (C1) actually provides a feasible region of $\zeta$
supported by the observations $D$, in which the proposed $\zeta$ has a
valid empirical likelihood value. If $Y_i < X_i^\top\beta(\tau_d)$ at
some $\tau_d$ for all $i=1,\ldots,n$, this proposed $\zeta$ will violate
the first condition, and then we regard its empirical likelihood value
as $0$. The second condition (C2) requires the set of estimating
functions to be linearly independent. Noting that
\[
E \{m(X,Y,\zeta_{0})m(X,Y,\zeta_{0})^\top\}= \Psi\otimes E(X
X^\top),
\]
where the elements of the $\Psi$ matrix\vspace*{2pt} are $\Psi_{ij} = \tau_i\wedge
\tau_j - \tau_i\tau_j$, the second condition is generally satisfied
for $\zeta$ near $\zeta_{0}$, as long as $E(X X^\top)$ is positive
definite.

For any proposed $\zeta$, consider its empirical likelihood function
$\mathcal{R}(\zeta)/n^n=\prod_{i=1}^n \omega_i(\zeta)$. With a prior
specification $p_0(\zeta)$ on the parameter $\zeta$, we can formally
have the posterior density
%e5 ###
%
\begin{equation}\label{eqa:postexp}
p (\zeta|D ) \propto p_0(\zeta) \times\mathcal{R}(\zeta) .
\end{equation}
We call $p (\zeta|D )$ the posterior distribution from the BEL
approach. This can be viewed as a misnomer, chosen for the sake of
convenience, because it is not really a posterior in the strict sense.
\citet{bellazar} proposed a procedure to check whether the empirical
likelihood is valid for posterior inference based on the criteria
provided in \citet{covproperty}. In this paper, we focus on the
asymptotic properties of the posterior distribution
(\ref{eqa:postexp}), and establish its frequentist validity by
first-order asymptotics.

Finding the maximum empirical likelihood estimator is a daunting task
computationally, because the objective function is generally
multi-modal. However, the value of the empirical likelihood ratio
$\mathcal{R}(\zeta)$ is relatively easy to compute given~$\zeta$, which
makes the Metropolis--Hastings algorithm, as given in
\citet{MHalgorithm}, feasible for sampling from the posterior. By
choosing a~proper prior, the posterior in (\ref{eqa:postexp}) is also
proper. Therefore, by checking the detailed balance equation and
Theorem 4.2 in \citet{MCMCconv}, the distribution of the MCMC sampler
converges to the posterior in (\ref{eqa:postexp}). More discussions on
computation efficiency can be referred to \citet{lteestMCMC}. A
Bayesian framework has its own merits in applications where informative
priors on $\beta(\tau)$ might be more appropriate than a strict
functional relationship on some of the parameters. For example, we may
believe that the slopes in $\beta(\tau_1 )$ are roughly the same as in
$\beta(\tau_{2})$. Imposing strict equalities to reduce the number of
unknown parameters in $\zeta$ might be hard to justify, but an
informative prior on the difference of two neighboring $\beta(\tau)$
can help regularize quantile estimation.

By using a standard Metropolis--Hastings algorithm for a given prior~$p_0(\zeta)$, we may use the average of the Markov chain on $\zeta$ as
an estimate of $\zeta$, when the posterior looks close to normal;
otherwise, we suggest using the mode of the posterior, which maximizes
(\ref{eqa:postexp}). In the empirical investigations in Sections
\ref{sec:simulation} and~\ref{sec:application}, we use the
posterior\vadjust{\goodbreak}
mode as the estimates.
%Large-sample inference on those parameters can be made in view of the
%asymptotic results obtained in the next section.
%The computational detail, which includes the evaluation of
%the empirical likelihood ratio function $\mathcal{R}(\zeta)$, choice
%of proposal distributions,
%and diagnostics of chain convergence, are not given in this paper, but
%are available in the supplementary material.

In our empirical investigations, we have found that the posterior mode
of the slope parameters behaves well, but the intercept parameter in
each~$\beta(\tau)$ can be better estimated in small samples if the
following strategy is followed. Suppose that
$\beta(\tau)=(\beta_I(\tau),\beta_S(\tau))$, where $\beta_I(\tau)$
corresponds to the intercept, and $\beta_S(\tau)$ corresponds to the
slope. Let ${\hat\beta}_S (\tau)$ be the posterior mode/mean obtained
from the MCMC chain, we use the modified estimate~${\hat\beta}_I
(\tau)$ as the $\tau$th sample quantile of $Y_i - X_{Si}^\top
{\hat\beta}_S (\tau)$, where $X_{Si}$ corresponds to $X_i$ excluding
the intercept term. This modification does not alter the asymptotic
distributions of the $\hat\beta(\tau)$. In the rest of the paper, we
always use this modification in the BEL estimate of quantile
regression.

%There are $k(p+1)$ estimating equations in (\ref{esteqael}) and
%$k(p+1)$ parameters in $\zeta$. If different quantiles share some
%common %features, we can reduce the dimensionality of $\zeta_k$. For
%example, if we assume the slope parameters across different quantiles
%to be the %same, which is true for regression models with
%homoscedastic errors, we can rewrite $\zeta$ as $(\beta_I(\tau_1),...,
%quantile $\tau$ and $\beta_C$ is the common slope across $k$
%quantiles. In this case, the dimension %of $\zeta$ is reduced to
%$k+p$. The empirical likelihood framework allows more estimating
%equations than the number of parameters. As shown %in \cite{empgee},
%the empirical likelihood method combines information from the
%estimating equations in an optimal way, that is, the empirical
%%likelihood method would automatically favor the optimal weights for
%linearly combining the estimation equations corresponding to different
%%quantiles.

%s3 ###
\section{Asymptotic properties of BEL} \label{sec:asymprop}

In this section, we provide an asymptotic justification of the BEL
estimator for quantile regression by deriving the limiting behavior of
the posterior distribution as $n \to\infty$. One noticeable point
about the estimating equations (\ref{esteqael}) is that they involve
indicator functions, so the resulting empirical likelihood ratio is
nonsmooth in $\zeta$. An asymptotic normality of the posterior
distribution in the Bayesian empirical likelihood context was derived
heuristically in \citet{bellazar} for smooth estimating equations. We
rely on empirical process theory to establish a similar result for the
BEL here.

As the first step, we shall prove the consistency of the maximum
empirical likelihood estimator (MELE), which is a necessary condition
for the asymptotic normality of the posterior.

%s3.1 ###
\subsection{Consistency of the MELE}
We assume that the true parameter $\zeta_0$ falls into a compact set of
the parameter space, and the optimization is carried out over this
compact set. For notational convenience, let
\[
\hat\zeta= \arg\max\{\mathcal{R}(\zeta)\}
\]
be the MELE, whose dependence on $n$ and the compact set on $\zeta$
have been suppressed in our notation. Note that the maximum empirical
likelihood estimate might not be unique, but the result here applies to
any maximizer of the empirical likelihood ratio, and all the maximizers
converge to the same asymptotic value.

The estimating functions $m(X,Y,\zeta)$ are not smooth in $\zeta$, but
it is worth noting that the expectations of $m(X,Y,\zeta)$ and the
empirical likelihood function are sufficiently smooth under the
following assumptions.

\begin{assumption}\label{aspA0}
There exists a neighborhood $\mathcal{N}$ of $\zeta_0$ such that\break
$P(\mathcal{R}(\zeta)>0 ) \rightarrow1$ for any $\zeta\in\mathcal{N}
$, as $n \to\infty$.
\end{assumption}
\begin{assumption}\label{aspA1}
The distribution function $G_X$ has bounded support~$\mathcal{X}$.
\end{assumption}
\begin{assumption}\label{aspA2}
The conditional distribution $F_X(t)$ of $Y$ given $X$ is twice
continuously differentiable in $t$ for all $X \in\mathcal{X}$.\vadjust{\goodbreak}
\end{assumption}
\begin{assumption}\label{aspA3}
At any $X \in\mathcal{X}$, the conditional density function\break
$F_X'(t)=f_X(t) > 0$ for $t$ in a neighborhood of $F_X^{-1}(\tau_d)$
for each $d =1, \ldots, k$.\vspace*{-1pt}
\end{assumption}

\begin{assumption}\label{aspA4}
$E\{m(X,Y,\zeta_{0})m(X,Y,\zeta_{0})^\top\}$ is positive definite.\vspace*{-1pt}
\end{assumption}

Assumption \ref{aspA0} is to guarantee that the interior of the convex
hull of $\{m(X_i,\allowbreak Y_i, \zeta)\dvtx i=1,\ldots,n\}$ for $\zeta\in\mathcal{N}$
contains the vector of zeros with probability tending to one. By
(\ref{exp-linearqrmodel}), $F_X(X^\top\beta_0(\tau_d))=\tau_d$ for any
$d \leq k$ and $X \in\mathcal{X}$. Therefore, for each $d$,
$\beta_0(\tau_d)$ is a solution to $E
\{m_{dk+j}(X,Y,\zeta)\}=0$, $j=0,\ldots,p$. Under Assumption
\ref{aspA3}, $\beta_0(\tau_d)$ is indeed the unique solution.
Correspondingly, $\zeta_{0}$ is the unique solution for $E
\{m(X,Y,\zeta)\}=0$.\vspace*{-1pt}

\begin{theorem}
\label{theorem-consistency} Under Assumptions \ref{aspA0}--\ref{aspA4}, the
MELE $\hat\zeta$ is a consistent estimator of $\zeta_{0}$.\vspace*{-1pt}
\end{theorem}

The proof of Theorem \ref{theorem-consistency} is sketched in the \hyperref[sec-proof]{Appendix}.
The basic idea is to check the conditions for consistency appearing in
Theorem 5.7 of \citet{asmstat}. Because those conditions require some
uniform convergence properties for collections of functions involving
$m(X,Y,\zeta)$, we use the empirical process theory as a natural tool.\vspace*{-1pt}

%s3.2 ###
\subsection{Asymptotic normality of the posterior}
To validate the asymptotic normality of the posterior distribution
(\ref{eqa:postexp}), we make one more assumption.\vspace*{-1pt}

\begin{assumption}\label{aspA5}
$\log\{p_0(\zeta)\}$ has bounded first derivative in a neighborhood of
$\zeta_0$.\vspace*{-1pt}
\end{assumption}

Then we have the following theorem.

\begin{theorem}\label{theorem-asynorm}
Under Assumptions \ref{aspA0}--\ref{aspA5}, the
posterior density of $\zeta$ has the following expansion on any
sequence of sets $\{\zeta\dvtx \zeta-\zeta_{0}=O(n^{-1/2}) \}$:
%e6 ###
%
\begin{equation}\label{postexpori}
p(\zeta|D) \propto\exp \bigl\{- \tfrac{1}{2}(\zeta-\hat\zeta
)^\top J_n(\zeta-\hat\zeta) + R_n \bigr\},
\end{equation}
where $\hat\zeta$ is the MELE,
\begin{eqnarray*}
J_n&=&nV_{12}^\top V_{11}^{-1}V_{12},\\
V_{11}&=& \Psi\otimes E(XX^\top),\\
V_{12}&=& - \frac{\partial E \{m(X,Y,\zeta) \}}{\partial\zeta}
\bigg|_{\zeta=\zeta_{0}}
\end{eqnarray*}
and $R_n=o_p(1)$. When $J_n$ is positive definite, we have
$J_n^{1/2}(\zeta-\hat\zeta)$ converging in distribution to
$N(0,I)$.
\end{theorem}

There are clear similarities between Theorem \ref{theorem-asynorm} here and
Theorem 1 of \citet{bellazar} for smooth estimating equations. We have
considered fixed priors, a common scenario in the literature, where the
limiting posterior\vadjust{\goodbreak} distributions of $\zeta$ are the same as the
limiting sampling distribution of the MELE [cf. \citet{empgee}]. An
important remark follows.\vspace*{-2pt}
\begin{remark}
The results in Theorem \ref{theorem-asynorm} apply to the cases where the
dimension of $\zeta$ is smaller than the dimension of the estimating
functions $m(X,Y,\zeta)$. For $\zeta$ with a reduced dimensionality,
the definition of $V_{12}$ is taken to be the derivative with respect
to the reduced parameter vector.\vspace*{-2pt}
\end{remark}

Asymptotically, Theorem \ref{theorem-asynorm} justifies the use of the BEL
approach\break for quantile regression with respect to frequentist
properties. When\break $f_X(X^\top\beta_0(\tau_d )) = f_{\tau_d}$ is constant
for all $X$, which is true for homoscedastic error models, we can
simplify $V_{12}$ to
\[
V_{12}= -\operatorname{diag} (f_{\tau_d})_{d=1,\ldots,k} \otimes E(XX^\top),
\]
if $\zeta$ is of $k(p+1)$ dimensions. Because $V_{11}=\Psi\otimes
E(\mathbf{X} \mathbf{X}^\top)$, the resultant asymptotic variance of
the posterior quantity, $J_n^{-1}$, is equivalent to the asymptotic
variance of the usual quantile regression (RQ) estimates, as proposed
in \citet{regquant}. This property is not shared by all working
likelihoods. If $\zeta$ is of lower dimensions, the posterior variance
no longer takes the same form, and improvements in the asymptotic
variances over RQ become possible.\vspace*{-2pt}

\begin{remark}
An improper prior cannot guarantee a proper posterior distribution. In
fact, the posterior will be improper for flat priors on $\zeta$ in the
BEL approach, and therefore we should avoid using flat priors on
$\zeta$.\vspace*{-2pt}
\end{remark}

Next, we consider a more interesting scenario where the prior
distribution shrinks with $n$. In this case, we use $p_{0,n}(\zeta)$ as
priors, and make the following assumption.\vspace*{-2pt}

\begin{assumption}\label{aspA6}
The logarithm of the prior density $p_{0,n}(\zeta)$ is twice
continuously differentiable, with the prior mode $\zeta_{0,n}=O(1)$,
and the matrix $J_{0,n} = -\frac{\partial^2 \log
\{p_{0,n}(\zeta)\}}{\partial\zeta^2} |_{\zeta=\zeta_{0,n}}=O(n)$.\vspace*{-2pt}
\end{assumption}

By Assumption \ref{aspA6}, $\log\{p_{0,n}(\zeta)\}$ can be Taylor
expanded up to the quad\-ratic term as follows.
%e7 ###
%
\begin{eqnarray}
\log\{p_{0,n}(\zeta)\}&=&\log\{p_{0n}(\zeta_{0,n})\}\nonumber\\ [-8pt]\\ [-8pt]
&&{} -
\tfrac{1}{2}(\zeta-\zeta_{0,n} )^\top
J_{0,n}(\zeta-\zeta_{0,n} )+ o (\Vert\zeta-\zeta_{0,n}\Vert^2).\nonumber
\end{eqnarray}
Then we have the following result.\vspace*{-2pt}

\begin{theorem}\label{theorem-asynorm2} Under Assumptions \ref{aspA0}--\ref{aspA4} and
\ref{aspA6}, the posterior density of $\zeta$ has the following
expansion on any sequence of sets $\{\zeta\dvtx
\Vert\zeta-\zeta_{0}\Vert=O(n^{-1/2}) \}$:
%e8 ###
%
\begin{equation}\label{postexpori2}
p(\zeta|D) \propto\exp \bigl\{- \tfrac{1}{2}(\zeta-\theta_{\mathrm{post}}
)^\top J_n(\zeta-\theta_{\mathrm{post}} ) + R_n \bigr\},\vadjust{\goodbreak}
\end{equation}
where
\begin{eqnarray*}
J_n&=&J_{0,n}+nV_{12}^\top V_{11}^{-1}V_{12},\\
\theta_{\mathrm{post}} &=& J_n^{-1} ( J_{0,n}\zeta_{0,n}+ n V_{12}^\top
V_{11}^{-1}V_{12}\hat\zeta)
\end{eqnarray*}
and $R_n=o_p(1)$.
\end{theorem}

Compared to Theorem \ref{theorem-asynorm}, the additional term $J_{0,n}$ in
both $J_n$ and $\theta_{\mathrm{post}}$ in Theorem \ref{theorem-asynorm2} provides a
balanced view of when and how an informative prior can complement the
likelihood in large samples. When $J_{0,n}=o_p(n)$, the posterior
expansion in Theorem \ref{theorem-asynorm2} is the same as that of Theorem
\ref{theorem-asynorm}, so the empirical likelihood will dominate the prior
information. Obviously, if~$J_{0,n}$ increases at a faster rate than
$n$, the prior will dominate the empirical likelihood. For the more
interesting case where $J_{0,n}$ increases at the rate of $n$, the BEL
produces a consistent estimate of $\zeta_0$ if $\Vert\zeta_{0,n} -
\zeta_0\Vert=o_p(1)$; otherwise, $\theta_{\mathrm{post}}$ may not converge to $
\zeta_0$ in probability, that is, a bias may be introduced, but the
variance is reduced. In the latter case, the posterior in
(\ref{postexpori2}) does not directly lead to asymptotically valid
posterior inference. However, noting that $J_n = J_{0,n} + nV_{12}^\top
V_{11} V_{12}$ and $J_{0,n}$ is known, the MCMC chain provides an
estimate of the matrix $nV_{12}^\top V_{11} V_{12}$, which is what we
need to obtain asymptotically valid confidence intervals.
%a consistent estimate of $\zeta_0$. Otherwise, bias will be
%introduced, but the variance of the BEL estimates will be smaller than
%that of the usual quantile regression estimates.
%yield asymptotically valid posterior inference. However, noting that
%$J_n = J_{0,n} + nV_{12}^\top V_{11} V_{12}$ and $J_{0,n}$ is known,
%the MCMC chain provides an estimate of the matrix $nV_{12}^\top V_{11}
%V_{12}$, which is what we need to obtain asymptotically valid
%confidence intervals.
%will dominate the empirical likelihood. If $||J_{0,n}||=o(n)$, the
%posterior expansion in Theorem \ref{theorem-asynorm} applies.
%%\item[C1]: If $\zeta_{0,n}$ does not shrink towards $\zeta_0$, the
%resultant $p(\zeta|D)$ may become a bi-modal distribution.
%%\item[C2]: If $J_{0,n}=o(n)$ instead of $O(n)$ in Assumption

Shrinking priors are relevant when the informative priors are
constructed from data of a secondary source or when the hypothesis on
common slope parameters are not rejected by a statistical test.

In Theorem \ref{theorem-asynorm2}, the prior mode $\zeta_{0,n}$ plays a
role in the posterior mean, which could be undesirable. For shrinking
toward common slopes, we can use a class of priors that eliminate the
bias due to a mis-specified prior mode when the common slope assumption
holds. For each $d=1, \ldots, k$, let $g_d$ be a spherically symmetric
distribution with zero as its center as well as its mode, and with a
finite second order derivative at zero. We consider a prior on $\zeta$
as
%e9 ###
%
\begin{eqnarray}\label{eqa-shrinkprcon}
\Omega^{-1/2}\bigl(\beta(\tau_1)-\beta_{p,0}\bigr) &\sim& g_1
\quad\mbox{and}\nonumber\\ [-8pt]\\ [-8pt]
\Sigma_d^{-1/2}\bigl(\beta(\tau_d)-\beta(\tau_1)\bigr) |
\beta(\tau_1) &\sim& g_d \qquad\mbox{for } d=2,\ldots,k \nonumber
\end{eqnarray}
for any location vector $\beta_{p,0}$ and scatter matrices $\Omega$ and
$ \Sigma_d$ of appropriate dimensions. They vary with $n$ in our
theory, but we have suppressed the dependence in notation. If we
write
\[
\Sigma_d =
\pmatrix{
\Sigma_{d, I} & \bzero^\top\cr
\bzero & \Sigma_{d, S}
}
,
\]
where $\Sigma_{d ,I} $ and $\Sigma_{d ,S}$ represent the components of
$\Sigma_d$ corresponding to the intercept and the slope parameters in
$\beta(\tau_d)$, respectively, for $d=2,\ldots,k$, we now assume
%e10 ###
%
\begin{equation}\label{eqa-priorvar}\qquad
\Vert\Omega^{-1}\Vert= O(\epsilon_n ),\qquad \Vert\Sigma_{d ,I}^{-1}\Vert =
O(\epsilon_n ) \quad\mbox{and}\quad \Vert\Sigma_{d,S} \Vert=O(n^{-1})
\end{equation}
for some sequence $\epsilon_n =o(n)$. We have the following
corollary.

\begin{corollary}\label{theorem-asynorm3}
Suppose that the same conditions of Theorem
\ref{theorem-asynorm2} hold. If the slope parameters in $\zeta_0$ are the
same at $\tau_1, \ldots, \tau_k$, and a (shrinking) prior satisfying
(\ref{eqa-shrinkprcon}) and (\ref{eqa-priorvar}) is used, the posterior
mean of Theorem \ref{theorem-asynorm2} becomes $\theta_{\mathrm{post}} = \zeta_0 +
O_p(\epsilon_n /n + n^{-1/2} ) .$
\end{corollary}

Clearly, Corollary \ref{theorem-asynorm3} indicates that the center of the
posterior is asymptotically unbiased for $\zeta_0$ with common slopes
regardless of what the prior mode $\beta_{p,0}$ is for $\beta(\tau_d )
$. All we need is to allow the prior variances of the slope differences
to be in the order of $1/n$, but the prior variances of the other
parameters increasing with $n$. The idea of constructing such a class
of shrinking priors applies more broadly than what we have considered
here with common slopes, but in our empirical work to be reported, only
independent normal and t-distributions will be used as $g_d$.

%s4 ###
\section{Simulation studies}\label{sec:simulation}

In this section, we use Monte Carlo simulations to investigate the
performance of the BEL methods (coverage probability and estimation
efficiency) from the frequentist viewpoint. We use the following
notation to distinguish BEL estimators with various priors on the
slope parameters. The usual quantile regression estimation at each
$\tau$ will be denoted simply as RQ.
\begin{itemize}
\item BEL.s: BEL estimators of single quantiles using moment
restrictions at each $\tau$.
%multiple $\tau$'s.
%assuming that the slope parameters vary linearly in $\tau$.
%
\item BEL.c: BEL estimators based on joint moment restrictions
assuming a~common slope parameter at several $\tau$'s.
\item BEL.n:
BEL estimators based on joint moment restrictions assuming that the
differences in slope parameters across $\tau$'s have normal priors with
zero mean and ``small'' variances.
\end{itemize}

%The common slope assumption and its variant BEL.n can be viewed as
%informative priors on the slope parameters. In the case of
%no commonality across quantiles, one can
%also compute empirical likelihood under joint moment restrictions at
%multiple $\tau$'s. In this case, the estimating functions at multiple
%quantiles use the same weights $\omega_i$, thus placing more
%restrictions on the feasible region for $\omega_i$. By Theorem
%estimators under single or joint moment restrictions are asymptotically
%equivalent to the RQ estimators.
%In finite samples, however, there is a noticeable difference between
%BEL.s and its counterparts based on joint
%moment restrictions. We have found that BEL.s generally performs
%better in empirical studies, and therefore, we have included only
%BEL.s when no assumptions are made on $\beta(\tau)$ across quantiles.

%In Section \ref{sec-coverage}, we examine the coverage probabilities
%and average lengths of the posterior intervals obtained by BEL.s. In
%%Section \ref{sec-efficiency}, we explore the estimation efficiency of
%these BEL estimators at tail quantiles $\tau=0.9,0.925,0.95$, to
%assess %the relative performances of the competing estimators. Noting
%that BEL.m did not show advantage over BEL.s in the simulation
%studies, we did %not report its performance in Section

%s4.1 ###
\subsection{Coverage properties}\label{sec-coverage} We first take a
brief look at the coverage probabilities of the
posterior credible intervals obtained under BEL.s. To see the impact of
empirical likelihood, we also include in the comparison two other
Bayesian methods, one based on the true parametric likelihood, and the
other based on a working likelihood.

%t1 ###
%
\begin{table}[b]
\caption{Comparison of $95\%$ posterior intervals of the median
regression parameters from three methods: \textup{(1)} BEL.s, \textup{(2)} BTL based on
the true likelihood and \textup{(3)} BDL based on a~working Laplace likelihood.
The coverage probability and lengths of the posterior intervals are
computed over $1000$ data sets of sample sizes $n=100, 400$ and
$1600$} \label{tab-coverage}
\begin{tabular*}{\textwidth}{@{\extracolsep{\fill}}lccccccc@{}}
\hline
& & \multicolumn{3}{c}{\textbf{Coverage of 95\% CI}} &  \multicolumn{3}{c@{}}{\textbf{Length of 95\% CI}}\\ [-5pt]
& & \multicolumn{3}{c}{\hrulefill} &  \multicolumn{3}{c@{}}{\hrulefill}\\
$\bolds{n}$ & & \textbf{BEL.s} &\textbf{BTL} &\textbf{BDL} & \textbf{BEL.s} &\textbf{BTL}
&\textbf{BDL}\\
\hline
\phantom{0}$100$ &$\beta_I(0.5)$ & 0.97 &0.94 &0.98 & 1.06 &0.80 &1.11 \\
&$\beta_S(0.5)$ & 0.98 &0.94 &0.98 & 0.58 &0.41 &0.58 \\
\phantom{0}$400$ &$\beta_I(0.5)$ & 0.97 &0.95 &0.98 & 0.43 &0.40 &0.55 \\
&$\beta_S(0.5)$ & 0.94 &0.95 &0.98 & 0.22 &0.20 &0.28 \\
$1600$ &$\beta_I(0.5)$ & 0.96 &0.96 &0.97 & 0.25 &0.21 &0.28 \\
&$\beta_S(0.5)$ & 0.96 &0.96 &0.98  &0.13 &0.10 &0.14 \\
\hline
\end{tabular*}
\end{table}

The data are generated from $Y_i = \beta_I + \beta_S (X_i -2) + e_i$
($i=1, \ldots, n $), where the true parameters are $\beta_I=2, \beta_S
=1$, $X_i$ and $ e_i $ are independently generated from the chi-square
distribution with 2 degrees of freedom and $N(0, 4)$, respectively. We
are interested in estimating the median regression coefficients
$\beta_I(0.5)$ and $\beta_S(0.5)$. Independent priors of $N(0, 100^2
)$, are used on both parameters.\vadjust{\goodbreak} We use the 2.5th and the 97.5th
percentiles of the Markov chain from BEL.s for $\tau=0.5$ to form 95\%
interval estimates for the parameters. The simulation study uses three
different sample sizes $n=100, 400, 1600$ to see whether the intervals
have desirable coverage probabilities for modestly large $n$.

In addition to BEL.s, we include two other Bayesian methods:
\begin{itemize}
\item BTL: the Bayesian method using the true likelihood
\[
\prod_{i=1}^n \sigma^{-1} \phi\biggl\{\frac{ y_i- \beta_I(0.5) -
\beta_S(0.5) (x_i-2) }{\sigma} \biggr\} ,
\]
where $\phi$ is the density of the standard normal distribution.
\item
BDL: a pseudo Bayesian method using the Laplace density as the working
likelihood
\[
\prod_{i=1}^n { \tilde\sigma}^{-1} \exp\biggl\{ -\frac{|y_i-
\beta_I(0.5) - \beta_S(0.5) (x_i-2)|}{2\tilde\sigma} \biggr\},
\]
where $\tilde\sigma$ is estimated by the mean of the absolute
residuals from the RQ estimate at $\tau=0.5$.
\end{itemize}

Similar MCMC sampling algorithms are used for all the three methods.
The BTL method can be viewed as a yardstick for any MCMC based method,
because it uses the true parametric likelihood under the model, which
is generally unknown in practice. The reason to consider BDL is that
the exponential component of its working likelihood is the objective
function of median regression. The BDL method has been used earlier by
\citet{bayesdeqr} among others, but in our empirical work, we have
chosen to use a fixed value of $\sigma$ in BDL, because we have found
that the MCMC chains have better mixing properties without including
$\sigma$ as an unknown parameter. A sensible value of $\sigma$ to use
in BDL is the RQ-based scale estimate. Table~\ref{tab-coverage}\vadjust{\goodbreak}
provides the average coverage probability and average length
information for each of the three methods over 1000 samples at each
choice of $n$.

% & &\multicolumn{3}{c}{\bf{$95\%$ CI}} & \multicolumn{3}{c}{\bf{LENGTH
%of $95\%$ CI}} \\
% \cline{3-8}
%{\bf{n}} & &\bf{BEL.s} &\bf{BTL} &\bf{BDL} &\bf{BEL.s} &\bf{BTL} &
%$100$ &a(0.5) &0.970 &0.944 &0.977 &1.057 &0.802 &1.106 \\
% &b(0.5) &0.976 &0.938 &0.982 &0.584 &0.409 &0.582 \\
%$400$ &a(0.5) &0.971 &0.953 &0.976 &0.425 &0.400 &0.554 \\
% &b(0.5) &0.943 &0.949 &0.981 &0.221 &0.199 &0.280 \\
%$1600$ &a(0.5) &0.956 &0.960 &0.973 &0.249 &0.213 &0.277 \\
% &b(0.5) &0.960 &0.956 &0.976 &0.126 &0.101 &0.139 \\
%regression parameters from three methods:
%(1) BEL.s, (2) BTL based on true likelihood, and (3) BDL based on a
%working Laplace likelihood. The coverage properties and lengths of are
%computed
%over 1000 data sets of sample sizes $n=100, 400, $ and 1600. }

This simple simulation study shows that as the sample size increases,
the posterior intervals obtained from BEL.s and BTL approach the
nominal levels $95\%$, although the convergence is not as fast as we
might have expected. Because the underlying model has i.i.d. normal
errors, the asymptotic relative efficiency of BEL.s and BDL are
approximately $67\%$ of BTL, which
helps explain the differences in the interval lengths. We also note
that BEL.s outperforms BDL by the frequentist measures, even after we
fixed the scale parameter in BDL.

Similar phenomena were observed in the interval estimation for other
quantiles and under several other error distributions, but we skip the
details. A more extensive report on estimation efficiency is given in
the next subsection.

%s4.2 ###
\subsection{Efficiency of BEL under various priors}\label{sec-efficiency}
In this section, we investigate the estimation efficiency of BEL.s,
BEL.c and BEL.n for $\zeta$ at different percentile levels, where the
posterior modes are taken as the parameter estimates. The estimation
efficiency is measured by the estimated mean squared error (MSE), with
data generated from the following four models:
\begin{itemize}
\item Model 1: $Y=X+Z+e$, where $X \sim\chi^2(2)$, $Z/2 \sim
\operatorname{Bernoulli}(0.5)$ and $e\sim N(0,4)$, with $X$, $Z$ and $e$ being
mutually independent;
\item Model 2: same as Model 1 except that
$\log(e) \sim N(0,1)$;
\item Model 3: $Y=X+Z+ (X/2+1)e$, where $X \sim
\chi^2(2)$, $Z/2 \sim \operatorname{Bernoulli}(0.5)$ and $e \sim N(0,4)$, with $X$,
$Z$ and $e$ being mutually independent;
\item Model 4: same as Model 3
except that $\log(e) \sim N(0,1)$.
\end{itemize}

These models include two covariates, of which $X$ is continuous, and
$Z$ is binary. Models 1 and 2 assume homoscedastic errors, and Models 3
and 4 allow the error distributions to depend on $X$. We use
$b_x(\tau)$, $b_z(\tau)$ to denote the two slope parameters, and
consider the adjusted intercept $a(\tau)$ as the fitted value of the
$\tau$th quantile at the sample mean of $(X,0)$. The reason that we
consider this adjusted intercept in the study, instead of the raw
intercept, is that the fitted value at the average design point of
$X_i$ is a more meaningful value than the fitted value at the origin,
which lies outside of the design space.
%a(\tau)=Q_{Y|\bar{X}}(\tau)=\bar{X}^\top\beta(\tau),

The three BEL methods (BEL.s, BEL.c and BEL.n) will be compared with RQ
and the composite quantile regression (CQR) of \citet{compositeqr}. The
CQR assumes common slopes, and minimizes the sum of individual quantile
loss functions over several $\tau$'s of interest. The CQR is a~direct
competitor of BEL.c, because they make the same assumption.

For Models 1 and 2, the common slope assumption holds, so there is no
asymptotic bias for any of the methods we consider here. Table
\ref{asymeff} shows the asymptotic efficiencies of BEL.c and CQR\vadjust{\goodbreak}
relative to RQ, when several quantiles are estimated jointly. It is
clear that BEL.c and CQR are similar in efficiency for Model~1, but
BEL.c stands out for Model 2. The asymptotic efficiency of BEL.s and
that of RQ are the same; both of them are improved on by the other
methods. Table \ref{asymeff} also includes comparisons at joint
estimation of three quartiles, to indicate that the efficiency gain of
BEL.c and CQR from the comparisons are not limited to high quantiles.

%t2 ###
%
\begin{table}
\caption{The table presents the ratio of the asymptotic MSE of the RQ
estimators over that of the BEL.c or CQR estimator for Models 1 and 2,
when jointly estimating quantiles at $\tau=0.25, 0.5, 0.75$ and
$\tau=0.9, 0.925, 0.95$, respectively}\label{asymeff}
\begin{tabular*}{\textwidth}{@{\extracolsep{\fill}}lcccccc@{}}
\hline
& \multicolumn{6}{c@{}}{\textbf{Asymptotic relative efficiencies for slope estimators}}
\\ [-5pt]
&\multicolumn{6}{c@{}}{\hrulefill}\\
 & $\bolds{\tau=0.25}$ & $\bolds{\tau=0.5}$ & $\bolds{\tau=0.75}$ &  $\bolds{\tau=0.9}$ & $\bolds{\tau
=0.925}$ & $\bolds{\tau=0.95}$ \\
\hline
 & \multicolumn{6}{c@{}}{Model 1} \\
%BEQ.l/RQ &0.959 &0.740 &0.959 &1.000 &0.883 &0.968\\
BEL.c$/$RQ & 1.598 &1.352 &\phantom{0}1.598 & 1.029 &1.219 &1.572\\
CQR$/$RQ & 1.590 &1.345 &\phantom{0}1.590 &  0.984 &1.166 &1.504\\ [6pt]
& \multicolumn{6}{c@{}}{Model 2} \\
%BEQ.l/RQ &0.996 &0.942 &0.606 & 0.999 &0.952 &0.883\\
BEL.c$/$RQ & 1.006 &3.280 &14.942 & 1.032 &1.677 &3.261 \\
CQR$/$RQ & 0.541 &1.763 &\phantom{0}8.032 & 0.756 &1.227 &2.386\\
\hline
\end{tabular*}\vspace*{-3pt}
\end{table}

The asymptotic efficiencies do not depend on the choices of fixed
priors. We now focus on estimation of high quantiles with
$\tau=0.9,0.925,0.95$ at the sample size of $n=100$, with the following
priors:
\begin{itemize}
\item For BEL.s and BEL.c, we use the prior $N(0,100^2)$ for each
intercept parameter, and $N(1,100^2)$ for each slope parameter.
\item
For BEL.n, we use the prior $N(0,100^2)$ for each intercept parameter,~and
$N(1,100^2)$ for $b_x(0.9)$ and for $b_z(0.9)$. The informative
priors used to~regulate the differences between quantiles are,
conditional on $\beta(0.90)$, $b_x(0.925) \sim N(b_x(0.9),0.16)$,
$b_x(0.95) \sim N(b_x(0.9),1)$, $b_z(0.925)\sim N(b_z(0.9),0.01)$ and
$b_z(0.95) \sim N(b_z(0.9) ,0.01)$.
\end{itemize}
Additional details of the Bayesian computations can be found in the
supplemental material [\citet{Yahe2012}]. The MSE's of various estimators of $\beta(\tau
)$ are given in Table \ref{mse-ext12} for Models~1 and 2, and in Table
\ref{mse-ext34} for Models~3 and 4. We make several observations from
those results:

\begin{itemize}
\item The performance of BEL.s is similar to or slightly better than
that of RQ.

\item When the common slope assumption holds, BEL.c has about the same
(Model~1) or better (Model 2) efficiency when compared with CQR. The
estimators that use informative priors on the slope parameters all
improve on RQ. The differences among various methods are more
significant at upper quantiles (say $\tau=0.95$) for heavier-tailed
distributions.
\item In Models 3 and 4, where the common slope
assumption does not hold for~$b_x (\tau)$, BEL.c\vadjust{\goodbreak} and CQR show
efficiency gains on the estimation of~$b_z (\tau)$, but losses in the
estimation of $b_x (\tau)$, due to bias. The BEL.n aims to reach a
compromise in the bias-variance trade-off, resulting in a better MSE
than RQ.
\end{itemize}

These findings are consistent with what we learned from the asymptotic
comparisons shown in Table \ref{asymeff}. The performance of BEL.n will
of course depend on the choice of priors on the difference in slopes.
The purpose of our study is not to demonstrate how to choose
informative priors, but to show how informative priors can make a
difference. Our empirical work shows that any reasonable choice of
priors helps, even though an optimal choice is too much to ask for in
general.\vspace*{-3pt}

%t3 ###
%
\begin{table}
\tabcolsep=0pt
\caption{The table gives the $n \times \mathit{MSE}$'s of several estimators for
the adjusted intercepts and slope parameters at three quantile levels
$\tau=0.9,0.925,0.95$ for Models 1 and 2, where $n=100$, and the $\mathit{MSE}$
is averaged over $500$ samples from each model. The numbers in the
brackets are the estimated standard errors}\label{mse-ext12}
\begin{tabular*}{\textwidth}{@{\extracolsep{\fill}}lccccccccc@{}}
\hline
& \multicolumn{3}{c}{\textbf{Adjusted intercepts}} & \multicolumn{6}{c@{}}{\textbf{Slopes}} \\ [-5pt]
& \multicolumn{3}{c}{\hrulefill} & \multicolumn{6}{c@{}}{\hrulefill} \\
\textbf{Method} & $\bolds{a(0.9)}$ & $\bolds{a(0.925)}$ & $\bolds{a(0.95)}$ & $\bolds{b_x(0.9)}$
& $\bolds{b_z(0.9)}$& $\bolds{b_x(0.925)}$ & $\bolds{b_z(0.925)}$ & $\bolds{b_x(0.95)}$ & $\bolds{b_z(0.95)}$\\
\hline
& \multicolumn{9}{c@{}}{Model 1}\\
BEL.s &22.0 &26.5 &35.6 &3.0 &11.7 &3.5 &13.8 &4.2 &19.8 \\
&(1.2) &(1.5) &(2.1) &(0.2) &(0.7) &(0.2) &(0.8) &(0.3) &(1.2) \\
%BEL.l &22.9 &26.7 &38 &3.6 &12.5 &3.6 &14.4 &4.9 &20.2 \\
%&(1.4) &(1.6) &(2.2) &(0.3) &(0.7) &(0.3) &(0.9) &(0.4) &(1.3) \\
BEL.n &23.1 &25.7 &31.5 &3.3 &12.3 &3.4 &12.3 &3.9 &12.4 \\
&(1.4) &(1.6) &(1.8) &(0.2) &(0.8) &(0.2) &(0.8) &(0.2) &(0.8) \\
BEL.c &26.6 &27.9 &34.1 &3.4 &13.9 &3.4 &13.9 &3.4 &13.9 \\
&(1.6) &(1.7) &(2) &(0.2) &(0.8) &(0.2) &(0.8) &(0.2) &(0.8) \\
CQR &22.8 &25.7 &30.0 &3.2 &12.7 &3.2 &12.7 &3.2 &12.7 \\
&(1.4) &(1.5) &(1.8) &(0.2) &(0.8) &(0.2) &(0.8) &(0.2) &(0.8) \\
RQ &22.3 &26.9 &36.5 &3.3 &12.1 &3.7 &14.4 &4.4 &19.2 \\
&(1.3) &(1.7) &(2.2) &(0.2) &(0.7) &(0.3) &(0.9) &(0.3) &(1.2) \\ [4pt]
& \multicolumn{9}{c@{}}{Model 2} \\
BEL.s &76.4 &126.6 &291.3 &9.5 &42.4 &13.5 &71.7 &26.2 &159.2\\
&(5.6) &(10.5) &(32) &(0.9) &(3.1) &(1.1) &(5.3) &(2.7) &(15.4) \\
%BEL.l &84 &139.1 &257.6 &8.5 &46.9 &12.9 &77.2 &23.8 &141 \\
%&(8.8) &(15.6) &(27.4) &(0.8) &(4.2) &(1.2) &(7.5) &(2.3) &(14.6) \\
BEL.n &78.7 &95.0 &150.0 &9.4 &43.6 &10.3 &43.8 &14.5 &43.7 \\
&(5.9) &(6.1) &(9.1) &(0.8) &(3.3) &(0.8) &(3.3) &(1.1) &(3.3) \\
BEL.c &86.8 &100.5 &158.3 &9.1 &46.9 &9.1 &46.9 &9.1 &46.9 \\
&(7.5) &(8.1) &(10.1) &(0.8) &(4.1) &(0.8) &(4.1) &(0.8) &(4.1) \\
CQR &109.3 &125.5 &175.9 &12.7 &61.7 &12.7 &61.7 &12.7 &61.7\\
&(11.1) &(11.3) &(15.5) &(1.2) &(5.2) &(1.2) &(5.2) &(1.2) &(5.2)\\
RQ &76.4 &136.4 &280.6 &10.0 &41.6 &14.9 &73.4 &26.5 &144.3\\
&(5.5) &(14.3) &(27.8) &(0.9) &(3.3) &(1.4) &(6.3) &(2.8) &(13.6)\\
\hline
\end{tabular*}\vspace*{-5pt}
\end{table}

%t4 ###
%
\begin{table}
\tabcolsep=0pt
\caption{Simulation results for Models 3 and 4; see the caption of
Table \protect\ref{mse-ext12} for more details}\label{mse-ext34}
\begin{tabular*}{\textwidth}{@{\extracolsep{\fill}}lccccccccc@{}}
\hline
& \multicolumn{3}{c}{\textbf{Adjusted intercepts}} & \multicolumn{6}{c@{}}{\textbf{Slopes}}
\\ [-5pt]
& \multicolumn{3}{c}{\hrulefill} & \multicolumn{6}{c@{}}{\hrulefill}\\
\textbf{Method} & $\bolds{a(0.9)}$ & $\bolds{a(0.925)}$ & $\bolds{a(0.95)}$ & $\bolds{b_x(0.9)}$ & $\bolds{b_z(0.9)}$
& $\bolds{b_x(0.925)}$ & $\bolds{b_z(0.925)}$& $\bolds{b_x(0.95)}$& $\bolds{b_z(0.95)}$\\
\hline
& \multicolumn{9}{c@{}}{Model 3} \\
BEL.s &90.2 &103.0 &138.7 &31.0 &35.9 &34.3 &42.6 &42.7 &66.2 \\
&(5.2) &(5.8) &(9.0) &(1.9) &(2.4) &(2.0) &(2.6) &(2.6) &(4.7) \\
%BEL.l &88.9 &96.6 &143.5 &36.9 &37.4 &36.1 &43.4 &47.4 &65.7 \\
%&(5) &(5.8) &(8.6) &(2.2) &(2.4) &(2.2) &(2.7) &(3) &(4.2) \\
BEL.n &95.6 &111.0 &129.4 &37.1 &41.4 &44.9 &41.5 &54.8 &41.7 \\
&(5.6) &(6.1) &(7.5) &(2.3) &(2.8) &(2.7) &(2.9) &(3.3) &(2.8) \\
BEL.c&104.3 &119.2 &143.0 &37.7 &43.5 &45.3 &43.5 &62.7 &43.5 \\
&(6.9) &(7.2) &(8.1) &(2.3) &(2.9) &(2.7) &(2.9) &(3.3) &(2.9) \\
CQR &94.6 &102.8 &118.0 &32.8 &38.4 &33.8 &38.4 &42.5 &38.4\\
&(5.9) &(6.3) &(7.3) &(2.0) &(2.6) &(1.9) &(2.6) &(2.4) &(2.6)\\
RQ &91.4 &106.9 &132.5 &30.6 &33.8 &35.0 &42.4 &42.9 &59.2\\
&(5.3) &(6.8) &(8.3) &(1.9) &(2.2) &(2.0) &(2.8) &(2.5) &(3.9)\\ [6pt]
& \multicolumn{9}{c@{}}{Model 4} \\
BEL.s &334.5 &507.5 &1085.1 &96.5 &134.1 &144.6 &213.5 &252.4 &547.6\\
&(25.7) &(40.0) &(109.5) &(8.4) &(10.1) &(12.3) &(17.6) &(19.9)
&(57.8)\\
%BEL.l &193.0 &280.3 &538 &85.9 &36.2 &105.8 &42.6 &174.4 &63.9 \\
%&(13.6) &(18.1) &(35.3) &(3.8) &(2.4) &(4.5) &(2.9) &(6.7) &(4.7) \\
BEL.n &277.0 &346.8 &518.0 &97.3 &124.7 &125.6 &124.7 &196.9 &125.1 \\
&(22.2) &(22.2) &(30.0) &(5.8) &(9.2) &(6.1) &(9.3) &(8.6) &(9.3)\\
BEL.c &391.6 &453.4 &659.8 &111.9 &160.5 &137.2 &160.5 &214.7 &160.5 \\
&(42.2) &(44.6) &(62.8) &(7.6) &(16.9) &(7.2) &(16.9) &(8.7) &(16.9) \\
CQR &530.1 &520.0 &663.9 &142.0 &195.3 &140.1 &195.3 &175.0 &195.3\\
&(56.4) &(52.1) &(58.2) &(12.0) &(18.8) &(10.3) &(18.8) &(9.5) &(18.8)\\
RQ &340.3 &552.0 &1014.6 &102.9 &123.9 &154.6 &215.0 &252.7 &481.5\\
&(25.9) &(56.6) &(101.1) &(8.3) &(9.6) &(11.5) &(19.0) &(17.8) &(46.0)\\
\hline
\end{tabular*}
\end{table}

%s5 ###
\section{An application to temperature downscaling} \label{sec:application}
In recent decades much focus has been placed on understanding potential
future climate changes. Meteorologists have developed various climate
models to simulate atmospheric variables\vadjust{\goodbreak} for both historical and future
time periods under different greenhouse gas emission scenarios.
Statistical downscaling approaches utilize those large-scale model
simulations to predict small-scale regional climate changes; see
\citet{statdown} for a review. Quantifying nearly extreme events in
climate studies is an important task, for which quantile regression is
a naturally appealing tool. However, high quantiles are usually hard to
estimate with RQ due to the inherently limited number of observations
in the tail of the distributions. In this section, we consider the BEL
methods for statistical downscaling of daily maximum temperature. We
used the observed daily maximum temperature (TMAX) of Aurora, IL
station from 1957--2002 as the response variable. The predictors are
the simulated daily maximum temperature (RTEM) and an indicator of wet
days (RAIN) from the ERA-40 reanalysis model introduced in
\citet{modeldata}. A~wet day is denoted by $\mathrm{RAIN} = 1$, when the
precipitation from ERA-40 is more than 1.2 kg/s/m$^2$. About 30\% of
the days are categorized as wet days in Aurora. We used the following
linear quantile regression model:
%e11 ###
%
\begin{equation}\label{equ:tempmodel}\qquad
Q_{\tau}(\mbox{TMAX} | \mbox{RTEM}, \mbox{RAIN}) = a(\tau) + b_x(\tau)
\mbox{RTEM} + b_z(\tau) \mbox{RAIN}\vadjust{\goodbreak}
\end{equation}
at high quantiles $\tau=0.99,0.995,0.999$. The quantile at $\tau=0.999$
is nearly extreme relative to our sample size, so the asymptotic theory
developed in this paper might be questioned. We choose to consider such
high quantiles partly to test the limits of our BEL methods.

We applied the following BEL methods with normal priors $N(0,1000^2)$
on each parameter to estimate the parameters of Model
(\ref{equ:tempmodel}), unless otherwise specified:
\begin{itemize}
\item BEL.c and BEL.s as introduced in Section \ref{sec:simulation}.
% \[
% \frac{\beta_S(0.999)-\beta_S(0.995)}{0.999-0.995} = \frac{
% \]
%$b_x(0.99)=b_x(0.995)=b_x(0.999)$.
%
\item BEL.z: the BEL estimator that assumes
$b_z(0.99)=b_z(0.995)=b_z(0.999)$.
%N(b_z(0.99),0.25^2)$, $b_x(0.999)|b_x(0.99) \sim N(b_x(0.99),0.6^2)$,
%and $b_z(0.999)|b_z(0.99) \sim N(b_z(0.99),2^2)$
%N(b_z(0.99),0.25^2)$, and $b_z(0.99)=b_z(0.995)=b_z(0.999)$.
%
\item BEL.t: the BEL estimator that assumes that given $b_x(0.99) $ and
$b_z(0.99)$, $(b_x(0.995)-b_x(0.99))/0.02$,
$(b_z(0.995)-b_z(0.99))/0.14$, $(b_x(0.999)-\break b_x(0.99))/0.35$ and
$(b_z(0.999)-b_z(0.99))/1.16$ are independent priors as the t
distribution with degrees of freedom $3$.
%t_3(b_z(0.99),0.35)$, and $b_z(0.99)=b_z(0.995)=b_z(0.999)$.
\end{itemize}
The scaling used in the prior distributions of BEL.t was chosen in
rough proportion to the variances of those parameter estimates from RQ,
and no optimality is claimed here. To assess the performances of
various methods, we randomly split the data from each year into two
parts, a fitting period and a testing period, with equal sizes of
$7889$ days in each part. We used the BEL methods and RQ for the
fitting period in estimating the model parameters and then applied the
fitted model to the testing period to predict the $\tau$th quantile of
TMAX. We randomly split the data three times, and labeled them as
SPLIT 1, SPLIT 2 and SPLIT 3, respectively. The average effective
sample sizes of the Markov chains for the BEL methods used here are
shown in Table \ref{tab:effectsz}, as calculated by the R function
effectiveSize() in the R package \textit{coda}.

%t5 ###
%
\begin{table}
\caption{Average effective sample sizes of the Markov chains used in
the downscaling example} \label{tab:effectsz}
\begin{tabular*}{\textwidth}{@{\extracolsep{\fill}}lccc@{}}
\hline
\textbf{Method} & \textbf{SPLIT 1} & \textbf{SPLIT 2} & \textbf{SPLIT 3} \\
\hline
BEL.c &976 &933 &364 \\
BEL.z &542 &633 &747 \\
BEL.t &672 &701 &515 \\
\hline
\end{tabular*}
\end{table}

%t6 ###
%
\begin{table}
\tabcolsep=0pt
\caption{The table presents the normalized differences calculated by
(\protect\ref{eqa:normdiff}). The row names provide the method used for model
fitting. In the column names, the Whole period indicates all the data
in the testing period are used; Lower RTEM indicates the testing data
with RTEM below its median; Wet~days indicates the testing data with
RAIN equals to $1$}\label{tab:predint}
\begin{tabular*}{\textwidth}{@{\extracolsep{\fill}}ld{2.3}d{2.3}d{2.3}d{2.3}d{2.3}d{2.3}d{2.3}d{2.3}d{2.3}@{}}
\hline
& \multicolumn{3}{c}{\textbf{Whole period}} & \multicolumn{3}{c}{\textbf{Lower RTEM}} & \multicolumn{3}{c@{}}{\textbf{Wet days}}
\\ [-5pt]
& \multicolumn{3}{c}{\hrulefill} & \multicolumn{3}{c}{\hrulefill} & \multicolumn{3}{c@{}}{\hrulefill}\\
\textbf{Method} & \multicolumn{1}{c}{$\bolds{>}$\textbf{0.99}} & \multicolumn{1}{c}{$\bolds{>}$\textbf{0.995}} &
\multicolumn{1}{c}{$\bolds{>}$\textbf{0.999}}
& \multicolumn{1}{c}{$\bolds{>}$\textbf{0.99}} & \multicolumn{1}{c}{$\bolds{>}$\textbf{0.995}}
& \multicolumn{1}{c}{$\bolds{>}$\textbf{0.999}}
& \multicolumn{1}{c}{$\bolds{>}$\textbf{0.99}} &
\multicolumn{1}{c}{$\bolds{>}$\textbf{0.995}} &
\multicolumn{1}{c}{$\bolds{>}$\textbf{0.999}}\\
\hline
& \multicolumn{9}{c@{}}{SPLIT 1} \\
%BEL.c &0.012 &0.089 &0.040 &-0.711 &-0.163 &1.036 &1.206 &0.425 &-0.316
BEL.c &0.012 &0.089 &0.040 &-0.871 &-0.163 &1.036 &0.402 &0.142 &-0.316
\\
BEL.z &0.012 &0.089 &0.040 &-0.230 &-0.163 &-0.476 &0.804 &0.142
&-0.316\\
%BEL.tp &0.012 &0.089 &0.040 &-0.871 &-0.388 &-1.483 &1.407 &0.709
%&0.949\\
BEL.t &0.012 &0.089 &0.040 & -1.351 &-0.163 &-1.483 &-0.201 & 0.142
&0.949\\
%BEL.cx&0.012 &0.089 &0.040 &-0.871 &-0.163 &0.532 &0.402 &0.142 &-0.949
%BEL.l &0.012 &0.089 &0.040 &-0.070 &-0.163 &0.532 &0.402 &0.142 &-0.949
%BEL.np &0.012 &0.089 &0.040 &-1.511 &-0.388 &-1.483 &0.804 &0.709
%&0.949 \\
%BEL.npx & 0.012 &0.089 &0.040 &-0.712 &-0.163 &-0.979 &1.005 &0.425
%&0.316 \\
%BEL.tpx &0.012 &0.089 &0.040 &-0.711 &-0.163 &-0.476 &1.005 &0.425
%&-0.316\\
%BEL.s &0.012 &0.089 &0.040 &-0.390 &-0.163 &-1.483 &0.804 &0.142 &0.949
RQ &\multicolumn{1}{c}{\textbf{$\bolds{-}$2.930}\phantom{\textbf{!}}} &\multicolumn{1}{c}{\textbf{$\bolds{-}$2.306}\phantom{\textbf{!}}}
&-1.742 &\multicolumn{1}{c}{\textbf{$\bolds{-}$2.151}\phantom{\textbf{!}}} &-1.743
&-1.987 &-1.005 &-1.276 &-1.582\\ [6pt]
& \multicolumn{9}{c@{}}{SPLIT 2} \\
%BEL.c &0.012 &0.089 &0.040 &-0.070 &0.515 &1.036 &\textbf{3.062}
%&1.875 &0.329 \\
BEL.c &0.012 &0.089 &0.040 &0.250 &0.515 &1.540 &\multicolumn{1}{c}{\phantom{!.}\textbf{2.659}} &1.591
&0.329 \\
BEL.z &0.012 &0.089 &0.040 &1.530 &-0.163 &0.532 &0.642 &0.452 &0.329 \\
%BEL.tp &0.012 &0.089 &0.040 &1.210 &1.192 &0.532 &1.448 &1.022 &0.329 \
BEL.t &0.012 &0.089 &0.040 & 1.050 &1.192 &0.532 &1.650 &1.022 &0.329 \\
%BEL.cx &0.012 &0.089 &0.040 &0.250 &0.740 &1.540 &1.852 &1.591 &0.329 \
%BEL.l &0.012 &0.089 &0.040 &1.210 &0.289 &0.532 &1.448 &1.022 &0.329 \\
%BEL.np &0.012 &0.089 &0.040 &1.210 &0.740 &0.532 &1.852 &1.591 &0.329 \
%BEL.npx &0.012 &0.089 &0.040 &0.730 &0.740 &0.532 &\textbf{2.054}
%&1.306 &0.329 \\
%BEL.tpx &0.012 &0.089 &0.040 &0.890 &0.966 &0.532 &1.852 &1.022 &0.329
%
%BEL.s &0.012 &0.089 &0.040 &1.370 &0.740 &0.532 &0.843 &1.306 &0.329 \\
RQ &0.012 &-0.390 &\multicolumn{1}{c}{\phantom{!.}\textbf{3.958}} &1.370 &0.063 &1.540 &0.843 &0.737
&\multicolumn{1}{c}{\phantom{!.}\textbf{4.139}} \\ [6pt]
& \multicolumn{9}{c@{}}{SPLIT 3} \\
%BEL.c &0.012 &0.089 &0.040 &-0.711 &0.063 &0.532 &-1.784 &-1.259
%&-0.308\\
BEL.c &0.012 &0.089 &0.040 &-1.511 &-0.163 &0.532
&\multicolumn{1}{c}{\textbf{$\bolds{-}$2.188}\phantom{\textbf{!}}}
&-1.543 &-0.308 \\
BEL.z &0.012 &0.089 &0.040 & 0.250 &-0.163 &-1.483 &-1.381 &-0.974
&0.962 \\
BEL.t &0.012 &0.089 &0.040 &-0.871 &-0.163 &-0.979 &-0.776 &-0.690
&1.596 \\
%BEL.cx &0.012 &0.089 &0.040 &-0.871 &-0.163 &0.532 &-1.986 &-1.259
%&-0.308 \\
%BEL.l &0.012 &0.089 &0.040 &1.370 &0.063 &0.532 &-0.776 &-1.259
%&-0.308 \\
%BEL.np &0.012 &0.089 &0.040 &-1.031 &-0.163 &-0.979 &-1.583 &-1.259
%&0.327 \\
%BEL.npx &0.012 &0.089 &0.040 &-0.871 &-0.163 &-0.980 &-1.179 &-0.974
%&0.962 \\
%BEL.tpx &0.012 &0.089 &0.040 &-0.711 &-0.163 &-1.483 &-1.784 &-1.259
%&0.962 \\
%
%BEL.s &0.012 &0.089 &0.040 &0.730 &-0.163 &-0.979 &-1.583 &-1.259
%&1.596 \\
RQ &-0.666 &-0.869 &\multicolumn{1}{c}{\textbf{$\bolds{-}$2.454}\phantom{\textbf{!}}} &0.250 &-0.388 &-1.483 &-1.583
&-1.259 &-1.577 \\
\hline
\end{tabular*}
\end{table}

Table \ref{tab:predint} reports the normalized differences as a
performance validation measure,
%e12 ###
%
\begin{equation}\label{eqa:normdiff}
d=\frac{O-E}{\sqrt{\tau(1-\tau)n}},
\end{equation}
where $n$ is the total number of days for prediction, $O$ is the number
of days when the observed TMAX exceeds the predicted $\tau$th quantile
of TMAX and $E$ indicates the expected number of days, that is,
$E=n(1-\tau)$. The normalized differences are shown for the whole
testing period, as well as for two subsets, one subset being the lower
half of RTEM, and the other subset being the wet days (RAIN~$=$ 1). The use
of these ad hoc subsets is meant to assess performances more
comprehensively. The normalized differences greater than $2$ in
absolute values are marked as bold in Table \ref{tab:predint}, from
which we have the following observations. First, over the whole testing
period, the normalized differences of each BEL method are stable across
random splits, but those from RQ predictions vary noticeably. For the
testing periods and for the selected subsets, the BEL methods perform
better than RQ, especially at $\tau=0.999$. Second, among the BEL
methods, BEL.c performs relatively worse, but BEL.t and BEL.z do well.
When we used the ANOVA test of \citet{anovatest} for the null
hypothesis of common slopes at $\tau=0.99, 0.995, 0.999$, the
hypothesis of $b_x(0.99)=b_x(0.995)=b_x(0.999)$ was rejected at 5\%
level of significance. This helps explain the inferior performance of
BEL.c relative to the other BEL methods, but all of them outperform RQ.

%We also fitted RQ directly to the testing data to obtain the fitted
%quantiles, and used them as a reference to measure the effectiveness
%of the predicted quantiles from the downscaling methods specified
%above.
%In reality, downscaling can give predictions to future TMAX, but the
%RQ fitted quantiles cannot be obtained for future time periods,
%so the study here may be viewed as a validation experiment.
%Table \ref{tab:predqua} reports the mean absolute difference between
%the predicted $\tau$-th quantiles from downscaling and the RQ fitted
%%quantiles
%to the testing data.
% We report the results for SPLIT 1 in the table because the results
%for SPLIT 2 and SPLIT 3 are similar.
%It is clear that as a downscaling method here, RQ predictions show
%higher deviations from the RQ fitted quantiles than the BEL methods,
%especially at $\tau=0.999$. Among the BEL methods, BEL.z has the
%smallest differences, followed by BEL.t, and then BEL.c.

Our empirical study shows that BEL methods can easily improve on RQ as
downscaling methods for high quantiles. Informative priors
will help further if the ``prior makers'' are well informed. In climate
studies, for example, historical data are generally available\vadjust{\goodbreak}
from multiple stations nearby, which can lead us to helpful
informative priors on slope parameters in the quantile models. In this
sense, the shrinking priors considered in Theorem \ref{theorem-asynorm2}
are relevant.

A natural question in climate downscaling is the autocorrelation of
measurements over time. In this section
we have bypassed this issue on two grounds. First, the quantile
regression estimation under the working
assumption of independence is typically consistent under weakly
dependent models; see \citet{longidep}. Second, we verified empirically
that the autocorrelation in TMAX was well represented by the
autocorrelations in the predictors used
in Model (\ref{equ:tempmodel}), and the signs of the residuals of the
quantile models were nearly uncorrelated. In more
general applications, however, it will be desirable to incorporate
dependence in an appropriate way, and future research
is clearly called for in this regard. Another interesting area of
future work is to perform downscaling at a group of stations and
include spatial correlation in the model. A recent paper by \citet
{ReichSpa} made a successful attempt at Bayesian spatial quantile
regression, and the idea of BEL with informative priors can be further
explored in spatial modeling.

%$\tau$-th quantile predictions from downscaling and the fitted
% regression quantiles on the testing data at $\tau=0.99, 0.995, 0.999$
%for SPLIT 1,2 and 3. }\label{tab:predqua}
% \bf{Method} & {\it\bf{$\tau=0.990$}} & {\it\bf{$\tau=0.995$}} & {
% \hline
% & \multicolumn{3}{c}{\textbf{SPLIT 1}} \\
% \cline{2-4}
%BEL.c &0.506 &0.283 &2.370 \\
%BEL.z &0.173 &0.142 &0.582 \\
%%BELtp & 0.618 &0.297 &1.202 \\
%BEL.t &0.591 &0.249 &2.035 \\

%BEL.cx &0.513 &0.165 &2.335 \\
%BEL.l &0.118 &0.143 &1.782 \\
%BEL.np &0.798 &0.311 &2.143 \\
%BEL.npx &0.462 &0.215 &0.226 \\
%BEL.tpx &0.453 &0.231 &0.280 \\
%BEL.s &0.264 &0.155 &1.743 \\
%RQ &0.821 &1.011 &2.539 \\
% & \multicolumn{3}{c}{\textbf{SPLIT 2}} \\
% \cline{2-4}
%BEL.c & 1.038 & 0.957 & 4.504 \\
%BEL.z & 0.792 &0.237 &1.606 \\
%BEL.p &0.953 &1.110 &2.883 \\
%BEL.t &1.002 &1.103 &3.157 \\
%BEL.cx & 0.603 & 0.901 & 4.654 \\
%BEL.l & 0.936 & 0.518 & 3.595 \\
%BEL.np & 1.091 &0.927 &3.719 \\
%BEL.npx &0.803 &0.860 &3.304 \\
%BEL.tpx &0.860 &0.947 &3.573 \\
%BEL.s & 0.823 & 0.704 & 2.184 \\
%RQ & 0.845 & 0.603 & 3.500 \\
% & \multicolumn{3}{c}{\textbf{SPLIT 3}} \\
% \cline{2-4}
%BEL.c & 0.596 &0.647 &2.765 \\
%BEL.c &0.947 &0.793 &2.402 \\
%BEL.z & 0.312 &0.457 &1.042 \\
%BEL.t & 0.472 &0.386 &0.797 \\
%BEL.cx & 0.687 & 0.584 & 2.809 \\
%BEL.l & 0.558 & 0.645 & 2.059 \\
%BEL.np & 0.588 & 0.558 & 0.413 \\
%BEL.npx & 0.479 & 0.44 & 0.809 \\
%BEL.tpx & 0.541 & 0.553 & 1.62 \\
%BEL.s & 0.408 & 0.614 & 0.7 \\
%RQ & 0.416 & 0.510 & 2.601 \\

%s6 ###
\section{Discussion} \label{sec:discussion}
In this paper, we propose using empirical likelihood as a~working
likelihood for quantile regression in Bayesian inference. We justify
the validity of the posterior based inference by establishing its first
order asymptotics. The BEL approach avoids the daunting task of
directly maximizing the EL function and allows informative priors to be
utilized. Although the idea of Bayesian quantile regression is not new,
the work provides an important addition to the literature by providing
the basic theory for incorporating possibly informative priors on
multiple quantiles. The efficiency gains are demonstrated through both
theoretical calculations and empirical investigations, when some common
features across quantiles are explored. If common slopes are assumed,
it is hard for the CQR method to find optimal weights in balancing the
quantile loss function at different $\tau$ levels, but the empirical
likelihood approach does so naturally. The use of informative priors is
also related in spirit to penalized optimization, but the lack of a
good overall objective function for several quantile levels makes the
usual regularization method difficult to formulate. The EL approach has
the ability to adapt automatically across quantile levels, and the BEL
approach enables flexible priors to be utilized in a simple way. Our
theoretical framework of shrinking priors provides good understanding
of how informative priors and likelihood can complement each other in
the BEL approach.

This paper uses empirical likelihood, but some of its variants such as
the ETEL, may work as well. The recent work of \citet{betelqr} provided
an approximation to the posterior from the Bayesian ETEL of quantile
regression at a given $\tau$. Although their approximation was not
strong enough to imply posterior\vadjust{\goodbreak} convergence for the Bayesian ETEL, it
can be strengthened using the approach we provide for BEL. We hope
that comparisons in a broader class of working likelihoods together
with efficient algorithms will be further developed in the
future.\vspace*{-3pt}

\begin{appendix}\label{sec-proof}
%s7 ###
\section*{Appendix: Proofs}
We begin with lemmas about the smoothness properties of functions
involving the estimating functions (\ref{esteqael}). Note that the
estimating functions~(\ref{esteqael}) involve an indicator function,
and as a result, the results obtained in \citet{empgee} for smooth
functions do not apply. While the work of \citet{empgee} relies on the
Taylor expansions, our proof uses the general theorem related to
M-estimators in \citet{asmstat} and the quadratic expansion
approximating the EL function provided in \citet{empnonsm}. We use
$x_j$ to indicate the $j$th component in the covarariates vector $X$
for $j=0,\ldots,p$, that is, $X=(x_0,x_1,\ldots,x_p)$ with $x_0=1$.

%s7.1 ###
\subsection{Preparatory results}
We discuss the properties of functions involving the estimating
function $m(X,Y,\zeta)$. Under Assumptions \ref{aspA1} and \ref{aspA2}
about~$G_X$ and $F_X$, $E \{m(X,Y,\zeta)\}$ can be sufficiently smooth.
\begin{lemma}\label{lemmaprop}
Under Assumptions \ref{aspA1} and \ref{aspA2}, we
have the following results:
\begin{longlist}
\item[(L1)] $E \{m(X,Y,\zeta)\}$ and
$E\{m(X,Y,\zeta)m(X,Y,\zeta)^\top \}$ are twice continuously
differentiable with respect to $\zeta$.
\item[(L2)] There exist
$k(p+1)$ dimensional compact neighborhoods $\mathcal{C_{\xi}}$ and
$\mathcal{C_{\zeta}}$ around $0$, in which
$E[m(X,Y,\zeta)/\{1+\xi^\top m(X,Y,\zeta)\}]$ is twice
continuously differentiable in $\zeta\in\mathcal{C_{\zeta}}$ and $\xi
\in\mathcal{C_{\lambda}}$, and
$E[m(X,Y,\zeta)m(X,Y,\zeta)^\top/\{1+\xi^\top m(X,\allowbreak Y,\zeta)\}]$
is uniformly continuous with respect to $\zeta\in\mathcal{C_{\zeta}}$
and $\xi\in\mathcal{C_{\lambda}}$.
\end{longlist}
\end{lemma}

\begin{pf}To show (L1), note that for each $d=0,\ldots,k-1$
and $j=0,\ldots,p$, there is
\begin{eqnarray*}
E\{m_{dk + j}(X,Y,\beta(\tau))\} &=& E \bigl\{\bigl(1_{\{ Y \leq X^\top\beta(\tau
_{d+1}) \} } - \tau_{d + 1} \bigr)x_j \bigr\} \\
&=& E_X \bigl[x_j \bigl\{E_{Y|X} \bigl(1_{\{ Y \leq X^\top\beta(\tau_{d+1}) \} } -
\tau_{d + 1} \bigr)\bigr\}\bigr] \\
&=& E_X [x_j \{F_X(X^\top\beta(\tau_{d+1}))-\tau_{d+1} \} ].
\end{eqnarray*}
Under Assumptions \ref{aspA1} and \ref{aspA2}, $E \{m(X,Y,\zeta)\}$ is
twice continuously differentiable. Consider the cases $i \leq l$ for
the second moments. By the definition of regression quantiles,
$X^\top\beta(\tau_i) \leq X^\top\beta(\tau_l)$, and therefore,
\begin{eqnarray*}
&&E \{m_{ik + j}(X,Y,\zeta)m_{lk + m}(X,Y,\zeta) \} \\
&&\qquad= E_X \bigl[x_j x_m \bigl\{E_{Y|X} \bigl(1_{\{ Y \leq X^\top\beta(\tau_{i + 1}
)\} }
- \tau_{i + 1} \bigr)\bigl(1_{\{ Y \leq X^\top\beta(\tau_{l + 1} )\} } - \tau
_{l + 1} \bigr)\bigr\} \bigr] \\
&&\qquad= E_X [x_j x_m \{F_X (X^\top\beta(\tau_{i + 1} )) - \tau_{l + 1}
F_X (X^\top\beta(\tau_{i + 1} ))\\
&&\qquad\hspace*{78pt}{}- \tau_{i + 1} F_X (X^\top\beta(\tau_{l + 1} )) + \tau_{i + 1}
\tau_{l + 1} \} ],
\end{eqnarray*}
which is twice continuously differentiable in $\zeta$.

Similarly, (L2) follows from
\begin{eqnarray*}
&&E \frac{m_{dk+j}(X,Y,\zeta)}{1+\xi^\top m(X,Y,\zeta)} \\
&&\qquad= E_X \biggl[
\sum_{0 \leq s \leq d} {\frac{{(1 - \tau_{d+1} )x_j }}
{{1 + \xi^\top m_s^* }}} \{F_X(X^\top\beta(\tau_{s + 1} )) -
F_X(X^\top\beta(\tau_s ))\} \\
&&\hspace*{53pt}{}- \sum_{d < s \leq k} {\frac{{\tau_{d+1} X_j }} {{1 + \xi
^\top m_s^* }}} \{F_X(X^\top\beta(\tau_{s + 1} )) -
F_X(X^\top\beta(\tau_s ) )\} \biggr],
\end{eqnarray*}
where we assume $\tau_0=0$, $\tau_{k+1}=1$, $m_0^*= ((1-\tau_1)
X^\top,\ldots,(1-\tau_k) X^\top)^\top$ and $m_s^* = ( -\tau_1
X^\top,\ldots,-\tau_s X^\top,(1 - \tau_{s + 1})X^\top,\ldots,(1 - \tau_k)
X^\top)^\top$ for $s=1,\ldots,k$. Because $m_s^*$ is bounded, $1+\xi
^\top m_s^*$ could be bounded away from $0$ for $\xi$ in a sufficiently
small compact neighborhood $\mathcal{C_{\xi}}$. Then $E
[m_{dk+j}(X,Y,\zeta_k)/\{1+\xi^\top m(X,Y,\zeta)\} ]$ is also twice
continuously differentiable in $\zeta$ and $\xi$. Similarly, we have
$E[m(X,Y,\zeta)m(X,Y,\zeta)^\top/\{1+\xi^\top m(X,Y,\zeta)\}]$
is uniformly continuous with respect to $\zeta\in\mathcal{C_{\zeta}}$
and $\xi\in\mathcal{C_{\lambda}}$.
\end{pf}

%
%s7.2 ###
\subsection{Consistency of the MELE}
By Assumptions \ref{aspA1}--\ref{aspA3}, the equation $E
\{m(X,Y,\zeta)\}=0$ has the unique solution $\zeta_{0}$. Define
%e13 ###
%
\begin{equation}\label{gamma}
\Gamma_n(\zeta)=
-n^{-1}\sum_{i=1}^{n}\log\{1+\lambda_n(\zeta)^\top
m(X_i,Y_i,\zeta)\},
\end{equation}
where $\lambda_n(\zeta)$ satisfies
\[
\sum_{i=1}^n\frac{m(X_i,Y_i,\zeta)}{1+\lambda_n(\zeta)^\top
m(X_i,Y_i,\zeta)}=0.
\]
Recall that
\[
\hat\zeta=\arg\max\{\Gamma_n(\zeta) \},
\]
we define the expected value of $\Gamma_n(\zeta)$ as
%e14 ###
%
\begin{equation}\label{expgamma}
\Gamma(\zeta)=-E [\log \{1+\xi(\zeta)^\top
m(X,Y,\zeta)\} ],
\end{equation}
where $\xi(\zeta)$ satisfies
\[
E
\biggl\{\frac{m(X,Y,\zeta)}{1+\xi(\zeta)^\top m(X,Y,\zeta)} \biggr\}=0.
\]
By Lemma \ref{lemmaprop}, Assumption \ref{aspA4}, and the implicit
function theorem, $\xi(\zeta)$ uniquely exists in the neighborhood
$\mathcal{C_{\lambda}}$ of $\underline{0} \in\mathbb{R}^{k(p+1)}$. To\vadjust{\goodbreak}
show that $\hat{\zeta}$ is a~consistent estimator of $\zeta_{0}$, it
is sufficient to check the conditions of Theorem~5.7 of
\citet{asmstat}. That is, we shall check
%e15 ###
%
\begin{equation}\label{cond1}
\sup_{\zeta} |\Gamma_n(\zeta)-\Gamma(\zeta)|
\stackrel{p}{\rightarrow} 0
\end{equation}
 and
\begin{equation}\label{cond2}
\sup_{\Vert\zeta-\zeta_{0}\Vert>\epsilon}
\Gamma(\zeta) < \Gamma(\zeta_{0})
\end{equation}
for any $\zeta$ within the compact neighborhood $\mathcal{C}_\zeta$ of
$\zeta_{0}$ and $\epsilon>0$.

\begin{lemma}
\label{concheck2} Under Assumptions \ref{aspA0}--\ref{aspA4},
(\ref{cond2}) holds.
\end{lemma}

\begin{pf}
It is easy to see $\xi(\zeta_0)=0$ because $E\{m(X,Y,\zeta_0)\}=0$, and
then $\Gamma(\zeta_0)=0$. By the Taylor expansion, we have
\[
\Gamma(\zeta)=-\xi(\zeta)^\top E
\biggl\{\frac{m(X,Y,\zeta)}{1+\xi(\zeta)^\top m(X,Y,\zeta)} \biggr\}
-\frac{1}{2}E\biggl\{\frac{(\xi(\zeta)^\top
m(X,Y,\zeta))^2}{(1+\alpha(\zeta)^\top
m(X,Y,\zeta))^2} \biggr\}
\]
for some $\alpha(\zeta)$ on the line segment between 0 and
$\xi(\zeta)$. On the right-hand side of the above equation, the first
term equals 0, and the second term with the negative sign included is
strictly negative, and thus $\Gamma(\zeta)<0$ for $\zeta\neq\zeta_0$.
So within the compact neighborhood $\mathcal{C}_\zeta$ of $\zeta_0$, we
have
\[
\sup_{\Vert\zeta-\zeta_0\Vert>\varepsilon}\Gamma(\zeta)<\Gamma(\zeta_0).
\]
\upqed\end{pf}

To check (\ref{cond1}), we first expand $\Gamma_n(\zeta)-\Gamma(\zeta)$
as
%e17 ###
%
\begin{eqnarray}\label{expansion}
\Gamma_n(\zeta)-\Gamma(\zeta) =Q_1 + Q_2,
\end{eqnarray}
where
\begin{eqnarray*}
Q_1&=& -n^{-1}\sum _{1 \leq i \leq n} [\log
\{1+\lambda_n(\zeta)^\top m(X_i,Y_i,\zeta)\} ]\\
&&{} + E [
\log\{1+\lambda_n(\zeta)^\top m(X_i,Y_i,\zeta)\} ],\\
Q_2 &=& - E [ \log\{1+\lambda_n(\zeta)^\top
m(X_i,Y_i,\zeta)\} ] + E [\log\{1+\xi(\zeta)^\top
m(X_i,Y_i,\zeta)\} ].
\end{eqnarray*}
To show the uniform convergence of (\ref{expansion}), we need the
following lemma.

\begin{lemma}\label{constvc} \textup{(i)} The class of constant functions:
$\mathcal{C}_0=\{\lambda,\lambda\in\mathcal{C}\}$ is
\mbox{P-Glivenko--Cantelli} (P-G--C) class, where $\mathcal{C}$ is some compact
set in $\mathbb{R}$. \textup{(ii)}~For bounded $X$, the class of functions
\begin{eqnarray*}
\mathcal{F}_1&=& \biggl\{ \frac{m(X,Y,\zeta)}{1+\lambda_n^\top
m(X,Y,\zeta)}\dvtx
\zeta\in\mathcal{C_\zeta},\lambda_n \in\mathcal{C_{\lambda}} \biggr\}\quad\mbox{and}\\
 \mathcal{F}_2&=&\bigl\{\log(\{1+\xi^\top
m(X,Y,\zeta)\}\dvtx \zeta\in\mathcal{C_\zeta}, \xi\in
\mathcal{C_{\lambda}}\bigr\}
\end{eqnarray*}
are P-G--C, where $C_{\lambda}$ is a compact neighborhood around
$\underline{0} \in\mathbb{R}^{k(p+1)}$, and $C_{\zeta}$ is a compact
neighborhood around $\zeta_0 \in\mathbb{R}^{k(p+1)}$.\vadjust{\goodbreak}
\end{lemma}

\begin{pf}
(i) According to Theorem 8.14 of \citet{emprocess} and the
fact that $\mathcal{C}_0$ is a collection of bounded functions, we only
need to show that $\mathcal{C}_0$ is VC-class, as defined in Section
9.1.1 in \citet{emprocess}. The P-measurability will be guaranteed by
the measurability and boundedness of the constant functions in $\mathcal
{C}_0$. The collection of all subgraphs of functions in $\mathcal{C}_0$
is $\mathcal{S}_0=\{(x,y),y<\lambda\}$. For any two points
$(x_1,y_1),(x_2,y_2)$ in $\mathbb{R}^2$, assume $y_1 \leq y_2$, it is
impossible that $\mathcal{S}_0$ would include $(x_2,y_2)$ while
excluding $(x_1,y_1)$. Therefore, based on the definition of
VC-subgraph Class, we have $\operatorname{VC}(\mathcal{C}_0)=2 < \infty$, i.e.,
$\mathcal{C}_0$ is a VC class. (ii) From Lemma~9.12 and Lemma 9.8 of
\citet{emprocess}, we know that the class of indicator functions
$\mathcal{G}_0 = \{ 1_{\{ Y \leq X^\top\beta\} },\beta\in{\mathbb
{R}}^{p+1} \}$ is a VC-class. From (vi) and (vii) in Lemma 9.9 of \citet
{emprocess}, the sets of estimating functions
\begin{eqnarray*}\label{eq:estfun}
\mathcal{G}_d = \bigl\{ \bigl(1_{\{ Y \leq X^\top
\beta(\tau_d) \} } - \tau_{i} \bigr)x_j, \beta(\tau_{d}) \in
{\mathbb{R}}^{p+1},0\leq j\leq p \bigr\},
\end{eqnarray*}
$1 \leq d \leq k$, are VC-class. Because $X$ is bounded,
$\mathcal{G}_d$ is P-G--C class by Theorem~8.14 of \citet{emprocess}.
Then by Theorem 9.26 of \citet{emprocess}, it follows that
$\mathcal{F}_1$ and $\mathcal{F}_2$ are P-G--C.
\end{pf}

We now verify (\ref{cond1}). We will check the uniform convergence of
$Q_1$ and~$Q_2$ in (\ref{expansion}).
Because $\mathcal{F}_2$, in which $\xi$ is not related to $(X,Y)$, is
P-G--C, the uniform convergence implied by P-G--C guarantees the
convergence of $Q_1$. For $Q_2$, because $\log\{1+\xi(\zeta)^\top
m(X_i,Y_i,\zeta))\}$ is bounded, by the dominate convergence theorem,
we only need to show $\lambda_n(\zeta)\stackrel{p}{\rightarrow} \xi
(\zeta)$ uniformly in $\zeta$. Because $\lambda_n(\zeta)$ is actually a
Z-estimator, the approximate zero of a data-dependent function of $\xi
(\zeta)$ as defined in Chapter 5.1 in \citet{asmstat}, then by using
the standard arguments of Z-estimator in \citet{asmstat} and by the
fact that $\mathcal{F}_1$ is P-G--C, we have $\lambda_n(\zeta) \stackrel
{p}{\rightarrow} \xi(\zeta)$ uniformly in $\zeta$.

The proof of Theorem \ref{theorem-consistency} is now complete.
%s7.3 ###
\subsection{Asymptotic normality of the posterior}
In our notation, we have
%e18 ###
%
\begin{equation}\label{empgamma}
\log\{ \mathcal{R}_n(\zeta) \} = n \Gamma_n(\zeta),
\end{equation}
where $\mathcal{R}_n(\zeta)$ is the empirical likelihood ratio of $\zeta
$. To expand $\Gamma_n(\zeta)$ up to the quadratic term, we use
Assumption \ref{aspA4}. We also use the following lemma, which is taken
from the quadratic expansion provided in Lemma A.6 of \citet
{empnonsm} %Molanes-L\'{o}pez, Van Keilegom and Veraverbeke (2008)
but formulated to suit our setting.

\begin{lemma}
\label{lemmaA6} Assume that the results of Lemma \ref{lemmaprop} and
Theorem \ref{theorem-consistency} hold. Under Assumptions
\ref{aspA0}--\ref{aspA4}, and additional conditions \textup{(C1)--(C3)} listed
below, we have
%e19 ###
%
\begin{eqnarray}\label{expgn}
\Gamma_n(\zeta) &=& -\tfrac{1}{2}(\zeta-\zeta_0)^\top V_{12}^\top
V_{11}^{-1}V_{12}(\zeta-\zeta_0) + n^{-1/2}(\zeta-\zeta_0)^\top
V_{12}^\top V_{11}^{-1}M_n \nonumber\hspace*{-35pt}\\ [-8pt]\\ [-8pt]
&&{}- \tfrac{1}{2}n^{-1}M_n^\top V_{11}^{-1}M_n + o_p(n^{-1})\nonumber\hspace*{-35pt}
\end{eqnarray}
uniformly in $\zeta$, for $\zeta-\zeta_0=O(n^{-1/2})$, and
%e20 ###
%
\begin{equation}\label{meleexp}
\hat\zeta-\zeta_0=n^{-1/2}(V_{12}^\top
V_{11}^{-1}V_{12})^{-1}V_{12}^\top V_{11}^{-1}M_n + o_p(n^{-1/2}),
\end{equation}
where $\hat\zeta$ is the MELE of $\zeta_0$, $M_n = n^{-1/2}
\sum_{i=1}^{n} m(X_i,Y_i,\zeta_0)$ and $V_{11}$ and $V_{12}$ are the
same as defined in Theorem \ref{theorem-asynorm}.
\end{lemma}

\begin{longlist}
\item[(C1)] $\Vert \sum_{i=1}^{n}
[m(X_i,Y_i,\zeta)-E\{m(X_i,Y_i,\zeta)\} ]\Vert=O_p(n^{1/2})$,
uniformly in $\zeta$ in a $o(1)$-neighborhood of $\zeta_0$.
\item[(C2)] $\Vert \sum_{i=1}^{n} [m
(X_i,Y_i,\zeta)m(X_i,Y_i,\zeta)^\top-E\{m(X_i,Y_i,\zeta)
m(X_i,Y_i,\zeta)^\top\} ] \Vert
= o_p(n)$, uniformly in $\zeta$ in a $o(1)$-neighborhood of $\zeta_0$.
\item[(C3)] $\Vert \sum_{i=1}^{n}
[m(X_i,Y_i,\zeta)-E\{m(X_i,Y_i,\zeta)\}\vspace*{1pt} -m(X, Y,\zeta_0)
+ E\{m(X,Y,\break \zeta_0)\}]\Vert
=o_p(n^{1/2})$, uniformly in $\zeta$ for $\zeta-\zeta_0=O_p(n^{-1/2})$.
\end{longlist}
%
%Noticing that $(C1)$ requires only the smoothness of the expectation
%of functions of $m(X,\mu,\nu)$.
To use the expansion (\ref{expgn}), we shall verify that (C1)--(C3)
are satisfied.

\begin{lemma}
\label{c456check} Under Assumptions \ref{aspA1}--\ref{aspA3},
Conditions \textup{(C1)--(C3)} are satisfied for the estimating functions
$m(X,Y,\zeta)$ of (\ref{esteqael}).
\end{lemma}

\begin{pf}
Because the collection of estimating functions $m(X,Y,\zeta)$ is\break
\mbox{P-Donsker} class, we have (C1). By the fact that the collection of the
product of the estimating functions is P-G--C, we have (C2). By
applying Lemma~$4.1$ of \citet{bahadur} to $m(X,Y,\zeta)$, we obtain
(C3).
%As shown in the proof of Lemma \ref{constvc}, the class $
%estimating functions $m_{ik+j}(X,Y,\zeta)$ in (\ref{eq:estfun}) is
%P-G--C. Then by Theorem 9.26 of \cite{emprocess}, the collection of the
%product of estimating functions is P-G--C. So $(C2)$ is also satisfied.
%By applying Lemma $4.1$ of \cite{bahadur} to $m_{ik+j}(X,Y,\zeta)$,
%$(C3)$ holds true.
\end{pf}

\begin{pf*}{Proof of Theorem \ref{theorem-asynorm}}
By Lemma \ref
{c456check}, Lemma \ref{lemmaA6}, (\ref{expgn}) and (\ref{empgamma}),
we have
\begin{eqnarray*}
\tilde p(\zeta|D) & = & p_0(\zeta) \times\mathcal{R}_n(\zeta) \\
&= & p_0(\zeta) \times\exp \biggl\{ -\frac{n}{2}(\zeta-\zeta_0)^\top
V_{12}^\top V_{11}^{-1}V_{12}(\zeta-\zeta_0)\biggr. \\
&&\biggl.\phantom{p_0(\zeta) \times\exp \biggl\{}{} +n^{1/2}(\zeta-\zeta_0)^\top V_{12}^\top V_{11}^{-1}M_n -
\frac{1}{2}M_n^\top V_{11}^{-1}M_n+ o_p(1) \biggr\}.
\end{eqnarray*}
Because of (\ref{meleexp}), we have
\begin{eqnarray*}
\tilde p(\zeta|D) &= & p_0(\zeta) \times\exp \biggl\{ -\frac{n}{2}(\zeta
-\zeta_0)^\top V_{12}^\top V_{11}^{-1}V_{12}(\zeta-\zeta_0)\biggr. \\
&&\biggl.\phantom{p_0(\zeta) \times\exp \biggl\{}{}+n(\zeta-\zeta_0)^\top V_{12}^\top V_{11}^{-1}V_{12}(\hat\zeta-\zeta
_0)\\
&&\hspace*{98pt}{} - \frac{1}{2}M_n^\top V_{11}^{-1}M_n + o_p(1) \biggr\} \\
& = & p_0(\zeta) \times\exp \biggl\{ -\frac{n}{2} (\zeta-\zeta_0)^\top
V_{12}^\top V_{11}^{-1}V_{12} (\zeta-2\hat\zeta+\zeta_0)\\
&&\hspace*{121pt}{} - \frac
{1}{2}M_n^\top V_{11}^{-1}M_n+o_p(1) \biggr\} \\
&= & p_0(\zeta) \times\exp \biggl\{ -\frac{n}{2} (\zeta-
\hat\zeta)^\top V_{12}^\top V_{11}^{-1}V_{12} (\zeta-
\hat\zeta) +o_p(1) \biggr\}.
\end{eqnarray*}

By Assumption \ref{aspA5}, we have
\[
\log\{p_0(\zeta)\} = \log\{p_0(\zeta_0)\} + O(n^{-1/2})
\]
for $\zeta-\zeta_0=O(n^{-1/2})$. Then we have
%e21 ###
%
\begin{equation}\label{postexpnorm2}
 \tilde p(\zeta|D) = p_{0}(\zeta_0) \exp \bigl\{
-\tfrac{1}{2} (\zeta- \hat\zeta)^\top J_n(\zeta- \hat\zeta)
+o_p(1) \bigr\},
\end{equation}
where $J_n= nV_{12}^\top V_{11}^{-1}V_{12}$.
For any $n$, we have $p(\zeta|D) \propto\tilde p(\zeta|D)$, and thus
(\ref{postexpori}) holds.

Because $J_n$ is positive definite, we have
%e22 ###
%
\begin{equation}\label{postnormexp}
p\bigl(J_n^{1/2}(\zeta-\hat\zeta)|D\bigr) \propto\exp
\bigl\{- \tfrac{1}{2} \bigl(J_n^{1/2}(\zeta-\hat\zeta)\bigr)^\top
\bigl(J_n^{1/2}(\zeta-\hat\zeta)\bigr) + o_p(1) \bigr\}\hspace*{-35pt}
\end{equation}
for any $\zeta-\zeta_0=O(n^{-1/2})$. Therefore, to show
\[
J_n^{1/2}(\zeta-\hat\zeta)\stackrel{D}{\rightarrow}N(0,I),
\]
it remains to show that
\[
P\bigl(\Vert J_n^{1/2}(\zeta-\hat\zeta)\Vert>\delta\bigr) \rightarrow0,
\]
when $\delta\rightarrow\infty$ and $n \rightarrow\infty$. From
(\ref{postexpnorm2}), we have for any $\zeta=\hat\zeta+ J_n^{-1/2}t$,
\[
\mathcal{R}_n(\zeta)\times p_0(\zeta) \stackrel{p}{\rightarrow}
p_0(\zeta_0)\exp\{-\Vert t\Vert^2/2\}.
\]

Because of $\mathcal{R}_n(\zeta) \times p_0(\zeta) \leq p_0(\zeta) $,
by the dominate convergence theorem, we have
\[
\int_{\Vert t\Vert >\delta} p_0(\hat\zeta
+J_n^{-1/2}t)\mathcal{R}_n(\hat\zeta+J_n^{-1/2}t)\,dt \rightarrow
p_0(\zeta_0)\int_{\Vert t\Vert > \delta}\exp\{-\Vert t\Vert^2/2\}\,dt
\]
for any $\delta\geq0$. Then it leads to
\begin{eqnarray*}
P\bigl(\Vert J_n^{1/2}(\zeta-\hat\zeta)\Vert >\delta| D\bigr) &=& \frac{\int_{\Vert t\Vert  >
\delta}
p_0(\hat\zeta+J_n^{-1/2}t)\mathcal{R}_n(\hat\zeta+J_n^{-1/2}t)\,dt}
{\int_{\Vert t\Vert  > 0} p_0(\hat\zeta+J_n^{-1/2}t)\mathcal{R}_n(\hat\zeta
+J_n^{-1/2}t)\,dt}\\
&\rightarrow& \frac{\int_{\Vert t\Vert  > \delta}\exp\{-\Vert t\Vert ^2/2\}\,dt}{\int
_{\Vert t\Vert  > 0}\exp\{-\Vert t\Vert ^2/2\}\,dt} \\
&=& (2\pi)^{-k(p+1)/2}\int_{\Vert t\Vert  > \delta}\exp\{-\Vert t\Vert ^2/2\}\,dt \\
&<& \epsilon
\end{eqnarray*}
for sufficiently large $\delta$.\vadjust{\goodbreak}
\end{pf*}

\begin{pf*}{Proof of Theorem \ref{theorem-asynorm2}}
Similar to the
proof of Theorem \ref{theorem-asynorm}, we have
%e23 ###
%
\begin{equation}\label{eqa:asynorm3}
\tilde p(\zeta|D) = p_{0,n}(\zeta) \times\exp \biggl\{ -\frac{n}{2}
(\zeta- \hat\zeta)^\top V_{12}^\top V_{11}^{-1}V_{12} (\zeta-
\hat\zeta) +o_p(1) \biggr\}.\hspace*{-35pt}
\end{equation}
By Assumption \ref{aspA6}, we have
\[
\log\{p_{0,n}(\zeta)\} = \log\{p_{0,n}(\zeta_{0,n})\} - \tfrac{1}{2}
(\zeta-\zeta_{0,n})^\top J_{0,n}(\zeta-\zeta_{0,n}) + o_p(1)
\]
for $\Vert\zeta-\zeta_0\Vert=O(n^{-1/2})$ and bounded $\zeta_{0,n}$. Combined
with (\ref{eqa:asynorm3}), we have
\[
\tilde p(\zeta|D) = C_n \exp \bigl\{-
\tfrac{1}{2}(\zeta-\theta_{\mathrm{post}} )^\top
J_n(\zeta-\theta_{\mathrm{post}} ) + R_n \bigr\},
\]
where $J_n=J_{0,n}+nV_{12}^\top V_{11}^{-1}V_{12}$, $\theta_{\mathrm{post}} =
J_n^{-1} ( J_{0,n}\zeta_{0,n}+ n V_{12}^\top
V_{11}^{-1}V_{12}\hat\zeta)$, $R_n=o_p(1)$, and $C_n$ is some
constant that does not depend on $\zeta$, and has the following
expression:
\[
C_n = p_{0,n}(\zeta_{0,n}) \exp \biggl\{-\frac{1}{2} \zeta_{0,n}^\top
J_{0,n}\zeta_{0,n} - \frac{n}{2}\hat\zeta^\top V_{12}^\top
V_{11}^{-1}V_{12}\hat\zeta+ \frac{1}{2} \theta_{\mathrm{post}}^\top J_n
\theta_{\mathrm{post}}\biggr \}.
\]
Therefore, we have (\ref{postexpori2}).
\end{pf*}
%
%It is easy to check the prior in (\ref{eqa-shrinkprcon}) is equivalent
%to
% \Omega&\Omega&\cdots&\Omega\\
% \Omega&\Sigma_2 + \Omega&\cdots&\Omega\\
% \vdots&\vdots&\ddots&\vdots\\
% \Omega&\Omega&\cdots&\Sigma_{k}+\Omega
% \end{pmatrix}.
%Due to the fact that
%J_{0,n} = \begin{pmatrix}
% \Omega^{-1}+ \Sigma_2^{-1}+...+\Sigma_k^{-1} &-\Sigma_2^{-1} &-
% -\Sigma_2^{-1} &\Sigma_2^{-1} &0 &\cdots&0 \\
% -\Sigma_3^{-1} &0 &\Sigma_3^{-1} &\cdots&0 \\
% \vdots&\vdots&\vdots&\ddots&\vdots\\
% -\Sigma_k^{-1} &0 &0 &\cdots&\Sigma_k^{-1}
% \end{pmatrix},
%$J_{0,n}$ increases at the rate of $n$, and therefore satisfies
%Assumption \ref{aspA6}. Noting that
%J_{0,n} (1_k \otimes\beta_{p,0}) = ( (\Omega^{-1}\beta_{p,0})^
%we have $||J_{0,n} (1_k \otimes\beta_{p,0}) ||=o_p(n)$. By Theorem

\begin{pf*}{Proof of Corollary \ref{theorem-asynorm3}}
The prior density $p_{0,n}(\zeta)$ can be written as
\begin{eqnarray*}
\log p_{0,n}(\zeta) &=& C + \log\bigl\{g_1\bigl(
\Omega^{-1/2}\bigl(\beta(\tau_1)-\beta_{p,0}\bigr)\bigr) \bigr\}\\
&&{} + \sum_{d=2}^k \log\bigl\{g_d\bigl(
\Sigma_d^{-1/2}\bigl(\beta(\tau_d)-\beta(\tau_1)\bigr)\bigr) \bigr\},
\end{eqnarray*}
where $C$ is some constant not depending on $\zeta$. Clearly, the prior
mode is $\beta(\tau_d)=\beta_{p,0}$ for all $d=1,\ldots,k$. Then we have
\[
\frac{\alpha^2 \log p_{0,n}(\zeta)}{\alpha\beta^2(\tau_1)}
\bigg|_{\zeta=1_k \otimes\beta_{p,0}} =\frac{\Omega^{-1/2}
g_1''(\bzero) \Omega^{-1/2}}{g_1(\bzero)} + \sum_{d=2}^k
\frac{\Sigma_d^{-1/2} g_d''(\bzero) \Sigma_d^{-1/2}}{g_d(\bzero)},
\]
and for $d=2,\ldots,k$,
\begin{eqnarray*}
\frac{\alpha^2 \log p_{0,n}(\zeta)}{\alpha\beta(\tau_1) \alpha
\beta(\tau_d)} \bigg|_{\zeta=1_k \otimes\beta_{p,0}}
&=&\frac{\alpha^2 \log\{g_d(\Sigma_d^{-1/2}(\beta(\tau_d) - \beta(\tau
_1)))\}}{\alpha\beta(\tau_1) \alpha\beta(\tau_d)} \bigg|_{\zeta=1_k
\otimes\beta_{p,0}} \\
&=& -\frac{\Sigma_d^{-1/2} g_d''(\bzero)
\Sigma_d^{-1/2}}{g_d(\bzero)},\\
\frac{\alpha^2 \log p_{0,n}(\zeta)}{\alpha\beta^2(\tau_d)}
\bigg|_{\zeta=1_k \otimes\beta_{p,0}}
&=&\frac{\alpha^2 \log\{g_d(\Sigma_d^{-1/2}(\beta(\tau_d) - \beta(\tau
_1)))\}}{\alpha\beta^2(\tau_d)} \bigg|_{\zeta=1_k \otimes\beta_{p,0}} \\
&=& \frac{\Sigma_d^{-1/2} g_d''(\bzero) \Sigma_d^{-1/2}}{g_d(\bzero)}.
\end{eqnarray*}
Note that for a spherically symmetric $g_d$ with its mode and center as
zero, we have
\[
\frac{g_d''(\bzero)}{g_d(\bzero)} = C_d I,
\]
where $I$ is the $(p+1)\times(p+1)$ dimensional identity matrix, and
$C_d >0 $ are constants for $d=1,\ldots,k$. Then, we have
\[
J_{0,n} =
\pmatrix{
\displaystyle C_1 \Omega^{-1} + \sum_{d=2}^k C_d \Sigma_d^{-1} &-C_2\Sigma_2^{-1}
&\cdots&C_k\Sigma_k^{-1} \cr
\displaystyle-C_2\Sigma_2^{-1} &C_2\Sigma_2^{-1} &\cdots&\bzero\cr
\vdots&\vdots&\ddots&\vdots\cr
\displaystyle-C_k\Sigma_k^{-1} &\bzero^\top &\cdots&C_k\Sigma_k^{-1}
}
\]
and therefore,
\begin{eqnarray*}
&&J_{0,n}(\zeta_{0,n} - \zeta_0)\\
&&\qquad =
\pmatrix{
\displaystyle C_1\Omega^{-1}\bigl(\beta_{p,0}-\beta_0(\tau_1)\bigr) + \sum_{d=2}^k C_d \Sigma
_{d,I}^{-1} \bigl(\beta_{0,I}(\tau_d)-\beta_{0,I}(\tau_1)\bigr)\cr
\displaystyle C_2\Sigma_{2,I}^{-1}\bigl(\beta_{0,I}(\tau_1)-\beta_{0,I}(\tau_2)\bigr)\cr
\vdots\cr
\displaystyle C_k\Sigma_{k,I}^{-1}\bigl(\beta_{0,I}(\tau_1)-\beta_{0,I}(\tau_k)\bigr)
}
,
\end{eqnarray*}
where $\beta_{0,I}(\tau_d)$ is the intercept parameter in
$\beta_0(\tau_d)$. Under the assumption in (\ref{eqa-shrinkprcon}) and
(\ref{eqa-priorvar}), $\Vert J_{0,n} (\zeta_{0,n} - \zeta_0)\Vert
=O(\epsilon_n)$ and $\Vert J_{0,n}\Vert $ is increasing at the rate of $n$.
Then Assumption \ref{aspA6} is satisfied, and Theorem
\ref{theorem-asynorm2} applies.

Note that the posterior mean $\theta_{\mathrm{post}}$ in Theorem
\ref{theorem-asynorm2} can be written as
\[
\theta_{\mathrm{post}} =\zeta_0 + n J_n^{-1} V_{12}^\top V_{11}^{-1} V_{12}
(\hat\zeta-\zeta_0) - J_n^{-1}J_{0,n}(\zeta_{0,n}-\zeta_0).
\]
By (\ref{meleexp}), we have $\Vert \hat\zeta-\zeta_0\Vert =O_p(n^{-1/2})$. Then
we have the posterior mean
\[
\theta_{\mathrm{post}} = \zeta_0 + O_p(\epsilon_n/n
+ n^{-1/2}).
\]
\upqed\end{pf*}
\end{appendix}

\section*{Acknowledgments}
The authors thank three anonymous referees and an
Associate Editor for their helpful and constructive comments and suggestions
that have led to several improvements in their work.

\begin{supplement}[id=suppA]
\stitle{Supplement to ``Bayesian empirical likelihood for quantile~regression''}
%suskaldyti doi
\sdatatype{.pdf}
\sfilename{aos1005\_supp.pdf}
\slink[doi]{10.1214/12-AOS1005SUPP}
\sdescription{The supplementary material contains additional details on the implementation
of the Bayesian computations used in the empirical studies reported in this paper.}
\end{supplement}

% imsref loaded by svajune.rapalyte, 2012-05-24 14:11:10
% imsref loaded by svajune.rapalyte, 2012-05-24 14:21:28
%

\printaddresses

\end{document}